\documentclass[12pt,reqno]{amsart}

\usepackage{amsthm,amsmath,amssymb,amsfonts}
\usepackage{epsfig}
\usepackage{psfrag,eepic}
\usepackage{url,hyperref}
\usepackage{ytableau}
\usepackage{tikz}
\usetikzlibrary{calc}

\makeatletter \@addtoreset{equation}{section}

\makeatother

\newcommand{\seteq}{\mathbin{:=}}

\newcommand{\Vir}{\mathrm{Vir}}

\newcommand{\Asl}{\widehat{\mathfrak{sl}}}
\newcommand{\rmF} {\mathrm{F}}
\newcommand{\rmQ} {\mathrm{Q}}
\newcommand{\rmV} {\mathrm{V}}
\newcommand{\bbL}{\mathbb{L}}
\newcommand{\bbF}{\mathbb{F}}
\newcommand{\calA}{\mathcal{A}}

\newcommand{\calF}{\mathcal{F}}

\newcommand{\calL}{\mathcal{L}}
\newcommand{\calM}{\mathcal{M}}
\newcommand{\calU}{\mathcal{U}}
\newcommand{\calV}{\mathcal{V}}

\renewcommand{\subsection}[1]{\refstepcounter{subsection}
\par\medskip\par\noindent\textbf{\thesubsection.~#1}}

\makeatletter
\def\Maketitle{{\def\newpage{}\maketitle}}
\def\Appendix{\appendix
  \def\@seccntformat##1{Appendix~\csname the##1\endcsname.~~}}
\makeatother

\newtheorem{Proposition}{Proposition}[section]

\newtheorem{Theorem}{Theorem}[section]

\newtheorem*{``Theorem''}{``Theorem''}

\theoremstyle{definition}
\newtheorem{Definition}{Definition}[section]

\newtheorem{Remark}{Remark}[section]

\begin{document}

\begin{title}[Coupling of two CFT and blow-up equations]
{Coupling of two conformal field theories and Nakajima-Yoshioka blow-up equations}
\end{title}

\author[M. Bershtein, B. Feigin and A. Litvinov]{Mikhail Bershtein, Boris Feigin and Alexei Litvinov}

\address{Landau Institute for Theoretical Physics, Chernogolovka, Russia
\newline 
Institute for Information Transmission Problems,  Moscow, Russia
\newline 
National Research University Higher School of Economics, Moscow, Russia
\newline 
Independent University of Moscow, Moscow, Russia}
\email{mbersht@gmail.com}

\address{National Research University Higher School of Economics, Moscow, Russia
\newline 
Landau Institute for Theoretical Physics, Chernogolovka, Russia
\newline 
Independent University of Moscow, Moscow, Russia}
\email{borfeigin@gmail.com}

\address{Landau Institute for Theoretical Physics, Chernogolovka, Russia
\newline 
NHETC, Department of Physics and Astronomy, Rutgers University, Piscataway, U.S.A.}
\email{litvinov@itp.ac.ru}



\begin{abstract}
We study the conformal vertex algebras which naturally arise in relation to the Nakajima-Yoshioka blow--up equations.
\end{abstract}

\Maketitle

\section{Introduction}

\subsection{}
Denote by $\mathit{M}(r,N)$ the moduli space of framed torsion free sheaves on $\mathbb{CP}^2$ of rank $r$, $c_1=0$, $c_2=N$. This space is a smooth partial compactification of the moduli space of $U(r)$ instantons. 

There is a natural action of the $r+2$ dimensional torus $
\mathit{T}$ on the $\mathit{M}(r,N)$: $(\mathbb{C^*})^2$ acts on the base $\mathbb{CP}^2$ and $(\mathbb{C^*})^r$ acts on the framing at the infinity. The Nekrasov partition function for the pure Yang-Mills theory is defined as the equivariant volume
\[Z(\epsilon_1,\epsilon_2,\vec{a};q)=\sum_{N=0}^{\infty} q^N \int_{\mathit{M}(r,N)}1,\]
where $\vec{a}=(a_1,\dots,a_r)$ and $\epsilon_1,\epsilon_2,a_1,\dots,a_r$ are the coordinates on the $\mathfrak{t}=\mathrm{Lie}\mathit{T}$. The last integrals can be computed by localisation method and equal to the sum of contributions of torus fixed points (which are labeled by $r$-tuple of Young diagrams $\lambda_1\dots,\lambda_r$). 

Nakajima and Yoshioka in paper \cite{NakajimaYoshioka} proved the so called blow--up equations for the function $Z(\epsilon_1,\epsilon_2,\vec{a};q)$. For the $r=1,2$ this equations have the form
\begin{align}
Z(\epsilon_1,\epsilon_2,a;q)&=Z(\epsilon_1,\epsilon_2{-}\epsilon_1,a;q)\cdot Z(\epsilon_1{-}\epsilon_2,\epsilon_2,a;q) \label{eq:blowr1} \\
	Z(\epsilon_1,\epsilon_2,a_1,a_2;q)&=\sum_{k \in \mathbb{Z}} \frac{q^{k^2}}{\mathrm{l}_k}Z(\epsilon_1,\epsilon_2{-}\epsilon_1,a_1{+}k\epsilon_1,a_2{-}k\epsilon_1;q)\cdot Z(\epsilon_1{-}\epsilon_2,\epsilon_2,a_1{+}k\epsilon_2,a_2{-}k\epsilon_2;q) \label{eq:blowr2}
\end{align}
The geometrical meaning of these equations is the relation between $\mathit{M}$ and $\widehat{\mathit{M}}$~--- the moduli space of framed torsion free sheaves on the blow--up of $\mathbb{CP}^2$. The shifted parameters $
(\epsilon_1,\epsilon_2{-}\epsilon_1
)$ and $(\epsilon_1{-}\epsilon_2,\epsilon_2)$ are weights of the torus action on the tangent space of two torus fixed points on the blow--up of $\mathbb{C}^2\subset \mathbb{CP}^2$.

For the $r=1$ the function $Z(\epsilon_1,\epsilon_2,a;q)=\exp(\dfrac{q}{\epsilon_1\epsilon_2})$ (see \cite[Sec 5]{NakajimaYoshioka}) so the equation \eqref{eq:blowr1} is trivial. But for $r=2$ the function $Z(\epsilon_1,\epsilon_2,a_1,a_2;q)$ coincide with the certain limit of the four point conformal block (due to the AGT conjecture \cite{AGT:2009}) for the conformal field theory (CFT) with the central charge $c=1+\dfrac{6(\epsilon_1+\epsilon_2)^2}{\epsilon_1\epsilon_2}$. Therefore the equation \eqref{eq:blowr2} suggests the relation in the conformal field theory. This relation is a main purpose of the paper.

\subsection{} We will mainly consider the $r=2$ case which corresponds to the CFT with the Virasoro algebra symmetry. Denote by $\mathrm{V}_{\Delta,c}$ the Verma module over the Virasoro algebra and by $\mathbb{L}_{\Delta,c}$ its irreducible quotient. The vacuum module $\mathbb{L}_{0,c}$ has the structure of the vacuum module of the Virasoro conformal vertex algebra. The corresponding operators are stress--energy tensor $T(z)$, its derivatives and products. We  parametrize central charge $c=c(b)=1+6(b+b^{-1})^2$ and denote this conformal vertex algebra as $\mathcal{M}_b$ similar to the notation for the minimal models. The representations $\mathbb{L}_{\Delta,c}$ become modules over the algebra $\mathcal{M}_b$. 

We know from the AGT conjecture that the left side of \eqref{eq:blowr2} is a conformal block\footnote{actually, the Whittaker-Gaiotto limit of the conformal block} in $\mathcal{M}_b$, where $b=\sqrt{\epsilon_1/\epsilon_2}$.. On the right side of \eqref{eq:blowr2} we have linear combination of the conformal blocks in $\calM_{b_1}\otimes \calM_{b_2}$, for the appropriate $b_1=b/\sqrt{b^2-1}$, $b_2=\sqrt{1-b^2}$.  It appears that the right side is a conformal block for some vertex algebra $\mathcal{A}_b$. This algebra $\mathcal{A}_b$ appears  to be an extension of the $\calM_{b_1}\otimes \calM_{b_2}$ by the field $\Phi_{1,3}\cdot \Phi_{3,1}$ (in notation of \cite{BPZ}). 

Geometrically two Virasoro algebras $\calM_{b_1}$ and $\calM_{b_2}$ correspond to two torus fixed points on the blow--up of $\mathbb{C}^2$. Algebraically one can extend the product $\calM_{b_1}\otimes \calM_{b_2}$ due to the relation between the central charges. In our case this relation is $b_1^2+b_2^{-2}=-1$, more general relations of this kind will be discussed in the Conclusion.

The identity \eqref{eq:blowr2} means the relation between the conformal vertex algebras $\mathcal{M}_b$ and $\mathcal{A}_b$. We prove for generic $b$ that $\mathcal{A}_b \cong\mathcal{M}_b\otimes \mathcal{U}$ for certain conformal vertex algebra $\mathcal{U}$, see Theorem \ref{th:Th1}.

This algebra $\calU$ is one of the main objects of the paper. It appears that $\calU$ can be constructed in terms of one free bosonic field $\varphi(z)$.  More precisely as a vertex algebra $\calU$ is isomorphic to lattice algebra $V_{\sqrt{2}\mathbb{Z}}$ or affine Lie algebra $\Asl(2)$ on the level 1. But $\calU$ has nonstandard (deformed) stress--energy tensor, its central charge equals to -5. see equation~\eqref{eq:Tu}. 

The algebra $\calU$ contains two commuting Virasoro subalgebras with central charges $\frac{-22}{5}$ and $\frac{-3}{5}$. This central charges correspond to minimal models $(2,5)$ and $(5,3)$. The sum of the corresponding stress--energy tensors $T_{2/5}+T_{3/5}$ equals  to the full stress--energy tensor $T_{\calU}$. As a consequence we can decompose $\calU$ as the tensor product of minimal models (see Theorem \ref{th:Ur} and character identities \eqref{eq:chi_sl2Vir2535}). 

\subsection{} The paper is organized as follows. In  Section \ref{Se:Res} we state the main results of the paper (which were shortly described above). The next section is devoted to proofs. The main tool of the proof is a Drinfeld--Sokolov reduction of the representation of $\Asl(2)_1\oplus\Asl(2)_k$ with respect to the diagonal $\Asl(2)_{k+1}$. Some arguments were based on the explicit computations which were made using Akira Fujitsu ope.math package \cite{ope.math}.

In the Section \ref{Seq:Comb} we discuss the combinatorial meaning of the character identities \eqref{eq:chi_sl2Vir2535} between  the characters of the $\Asl(2)$ of the level 1 and minimal models $(2,5)$ and $(5,3)$. In the Section \ref{Se:Equations} we show how to use our results to the blow--up equations. 

In the Conclusion we discuss possible generalizations of the product $\calM_{b_1}\otimes \calM_{b_2}$ and the algebra $\calU$. 

%

\section{Results} \label{Se:Res}
\subsection{Vertex algebras.}  We will use the language of vertex algebras, see e.g. \cite{FrenkelBook}. Recall that a vector space $V$ is called a vacuum representation of vertex algebra if any vector $v\in V$ corresponds to a power series of operators $Y(v;z)=\sum Y_nz^{-n}$, $Y\in \mathrm{End}(V)$. This correspondence $v \leftrightarrow Y(v;z)$ is called the \emph{operator-state correspondence}. In the definition of the vertex algebra the correspondence $v \longleftrightarrow Y(v;z)$ should satisfy certain conditions: vacuum axiom, translation axiom and locality axiom.

Recall that the vertex algebra $V$ is called \emph{conformal} if there exists a non-zero \emph{conformal vector} $\omega\in V$ such that corresponding power series of operators $T(z)$ satisfy
\begin{align} \label{eq:T(z)OPE}
T(z)T(w)=\frac{c}{(z-w)^4}+\frac{2}{(z-w)^2}T(w)+\frac{1}{(z-w)}\partial T(w)+\mathrm{reg}.
\end{align}
The corresponding $T(z)$ is called the \emph{stress--energy tensor}, parameter $c$ is called the \emph{central charge}. If we expand $T(z)$ into power series $T(z)=\sum_n L_n z^{-n-2}$ then equation \eqref{eq:T(z)OPE} is equivalent to the Virasoro algebra ($\Vir$ for short) relations 
\[[L_n,L_m]=(n-m)L_{n+m}+\frac{n^3-n}{12}c \delta_{n,-m}.\]

\subsection{Algebra $\mathcal{U}$.} \label{subU}
First we recall the construction of the lattice vertex algebra for one-dimensional lattice $\sqrt{2}\cdot \mathbb{Z}$ (see \cite[Sec 5.2]{FrenkelBook}). Let $a_n$ be the generators of the Heisenberg algebra
\[[a_n,a_m]=n\delta_{m+n,0}.\]
It is convenient to consider the operators $a_n$ as modes of 	the bosonic field $\varphi(z)$
\begin{align} \label{eq:varphi}
\varphi(z)=\sum_{n \in \mathbb{Z}\setminus 0} \dfrac{a_n}{-n} z^{-n}+a_0\log z+\widehat{Q}, 
\end{align}
where the operator $\widehat{Q}$ is conjugate to the operator $\widehat{P}=a_0$, i.e. satisfy the relation $[\widehat{P},\widehat{Q}]=1$.
The relations of the Heisenberg algebra can be rewritten in terms of the operator product expansion \[ \varphi(z)\varphi(w)= \log(z-w)+\mathrm{reg}.\] We will contract such notation to 
$\varphi(z)\varphi(w)\sim \log(z-w)$ below. 

Denote by $\rmF_{\lambda}$ the Fock representation of the Heisenberg algebra with the highest weight vector $v_\lambda$
\begin{align*}
a_nv_\lambda=0\;\, \text{for $n>0$},\quad a_0v_\lambda=\lambda v_\lambda.
\end{align*}
Denote by $S_\lambda$ the shift operator $S_\lambda\colon \rmF_{\mu} \rightarrow \rmF_{\mu+\lambda}$ defined by 
\[ S_\lambda v_{\mu}=v_{\mu+\lambda},\qquad [S_\lambda,a_n]=0,\text{ for } n\neq 0
\]
Actually $S_\lambda$ is just an exponent $\exp(\lambda \widehat{Q})$.

The direct sum $V_{\sqrt{2}\mathbb{Z}}\seteq \bigoplus_{k\in \mathbb{Z}}\rmF_{k\sqrt{2}}$ has a vertex algebra structure. This algebra is called \emph{the lattice vertex algebra} for the lattice $\sqrt{2}\cdot\mathbb{Z}$. Under the operator-state correspondence the highest weight vectors $v_\lambda$, $\lambda=k\sqrt{2}$ correspond to
\[Y(v_{\lambda};z)=:\!e^{\lambda\varphi}:=S_{\lambda}z^{\lambda a_0} \exp\Bigl(\lambda \sum_{n \in \mathbb{Z}_{>0}} \dfrac{a_{-n}}{n} z^{n} \Bigr) \exp\Bigl(\lambda\sum_{n \in \mathbb{Z}_{>0}} \dfrac{a_n}{-n} z^{-n} \Bigr), \]
Here and below $:\ldots:$ denotes the creation-annihilation normal ordering. 
For more general vectors of the form $v=a_{-m}^{n_m}\cdots a_{-1}^{n_1}v_{\lambda}$ the corresponding operators have the form
\[Y(v;z)=:\!(\partial^m\varphi)^{n_m} \cdots(\partial\varphi)^{n_1} e^{\lambda\varphi} :  \]
For more details about this construction see \cite[Sec. 5.2]{FrenkelBook}. 

The  algebra $V_{\sqrt{2}\mathbb{Z}}$ is isomorphic to the vertex algebra of affine Lie algebra $\Asl(2)$ on the level 1. We denote the standard generators of  $\Asl(2)=\mathfrak{sl}(2)\otimes \mathbb{C}[t,t^{-1}]\oplus \mathbb{C}K$ by $e_n=e\otimes t^n$, $f_n=f\otimes t^n$, $h_n=h\otimes t^n$, and the central element by $K$. 
We denote by $\mathcal{L}_{h,k}$ the irreducible module of  the $\widehat{\mathfrak{sl}}(2)$ algebra generated by the highest vector $v$ such that
\begin{align*}
e_nv=0,\,\, \text{for $n\geq0$};\qquad f_nv=h_nv=0,\,\, \text{for $n>0$};\qquad h_0v=hv,\,\, Kv=kv.
\end{align*}
The value $k$ of the central element is called the level of the representation.

The module $\mathcal{L}_{0,k}$ has the structure of the vacuum module of the vertex algebra. This algebra for generic $k$ will be used in the next section. If $k=1$ then the vertex algebra $\mathcal{L}_{0,1}$ is isomorphic to the lattice algebra $V_{\sqrt{2}\mathbb{Z}}$. The action of the generators of $\Asl(2)$ is defined by the formulas
\[
\sum_{n \in \mathbb{Z}}\! e_nz^{-n-1}=:\!e^{\sqrt{2}\varphi}:,\qquad
\sum_{n \in \mathbb{Z}}\! f_nz^{-n-1}=:\!e^{-\sqrt{2}\varphi}:,\qquad \sum_{n \in \mathbb{Z}}\! h_nz^{-n-1}=\sqrt{2} \partial \varphi(z).
\]

The standard conformal vector for the algebra $V_{\sqrt{2}\mathbb{Z}}=\mathcal{L}_{0,1}$ is $\omega_0=\frac12a_{-1}^2v_0$, the corresponding stress--energy tensor equals $T(z)=\frac12(\partial\varphi)^2$, and has the  central charge $1$. This vertex algebra 
has two representations: the vacuum representation and the second one $\bigoplus_{k\in \mathbb{Z}+
1/2}\rmF_{k\sqrt{2}}=\mathcal{L}_{1,1}$. Their characters i.e. traces of $q^{L_0}$ are
\begin{align*}
 \chi(\mathcal{L}_{0,1})=\chi\left(\bigoplus_{k\in \mathbb{Z}}\rmF_{k\sqrt{2}}\right)= \sum_{k \in \mathbb{Z}} {q^{k^2}}/{(q)_\infty}=1+3q+4q^2+7q^4+\dots,\\
 \chi(\mathcal{L}_{1,1})=\chi\left(\bigoplus_{k\in \mathbb{Z}+\frac12}\rmF_{k\sqrt{2}}\right)= \sum_{k \in \mathbb{Z}+\frac12} {q^{k^2}}/{(q)_\infty}=2q^{1/4}+2q^{5/4}+6q^{9/4}+\dots, 
\end{align*}
 where we used that $L_0v_\lambda=(\lambda^2/2)v_\lambda$ and $(q)_\infty=\prod_{k=1}^\infty (1-q^k)$.
 
There are other possible conformal vectors in $V_{\sqrt{2}\mathbb{Z}}$. Namely the local operators $\frac12(\partial\varphi)^2 +u(\partial^2\varphi)$ satisfy stress--energy relation \eqref{eq:T(z)OPE} with the central charge $c=1-12u^2$. The addition of $u(\partial^2\varphi)$ changes the $L_0$ operator to $L_0-u a_0$. In particular the eigenvalues of the new $L_0$ are integers if and only if $u \in \frac1{\sqrt{2}}\mathbb{Z}$.
 
Now we can define $\calU$.
\begin{Definition} The conformal vertex algebra $\mathcal{U}$ coincides with the $V_{\sqrt{2}\mathbb{Z}}$ as the vertex algebra, but the stress--energy tensor is modified
\begin{multline}
T_{\mathcal{U}}=\frac12(\partial\varphi)^2 +\frac1{\sqrt{2}}(\partial^2\varphi) +\epsilon\left(2(\partial\varphi)^2e^{\sqrt{2}\varphi} +\sqrt{2}(\partial^2\varphi)e^{\sqrt{2}\varphi}\right)= \\
=\frac12\partial_z\varphi(z)^2+\frac1{\sqrt{2}}\partial_z^2\varphi(z) +\epsilon \partial_z^2e(z), \quad \varepsilon\neq 0 \label{eq:Tu}
\end{multline}
\end{Definition}

It is clear that the conformal vertex algebras $\calU$ are isomorphic for different values $\varepsilon \neq 0$. For $\varepsilon=0$  $T_\calU(z)$ has the form discussed above  for $u=\frac{1}{\sqrt{2}}$ and central charge $-5$. 

The additional term in \eqref{eq:Tu} corresponds to the vector  $2e_{-3}v_0=(2a_{-1}^2+\sqrt{2}a_{-2})v_{\sqrt{2}}$.  The eigenvalue of the operator $L_0-\frac{1}{\sqrt{2}}a_0$ on this vector equals to 2. The OPE of the~corresponding operator $\partial_z^2 e(z)$ has no singular terms (since $e(z)e(w)\sim (z-w)^2$). Therefore $T_{\calU}$ defined in \eqref{eq:Tu} satisfies stress--energy tensor OPE \eqref{eq:T(z)OPE} with the central charge $c_{\calU}=-5$.

We will call this conformal vertex algebra the Urod algebra. This algebra is one of the main objects of the paper.

The representations of $\calU$ are the same as the representations   of $V_{\sqrt{2}\mathbb{Z}}$ $U_0=\bigoplus_{k\in \mathbb{Z}}\rmF_{k\sqrt{2}}$ and $U_1=\bigoplus_{k\in \mathbb{Z}+1/2}\rmF_{k\sqrt{2}}$, but with the different characters due to shift of the $L_0$ operator. Schematically this shift of the grading is represented on the following picture

\begin{figure}[h]
\begin{center} \small
\begin{tikzpicture}
\draw (-1,0) node {$v_0$};
\draw (1,-1) node {$v_{\sqrt{2}}$} edge[bend right,->,densely dotted] (1,0) ;
\draw (-1,-1) node {$a_{-\!1}v_{0}$};
\draw (-3,-1) node {$v_{-\!\sqrt{2}}$} edge[bend right,->,densely dotted] (-3.2,-1.8) ;
\draw (1,-2) node {\scriptsize $a_{-\!1}v_{\sqrt{2}}$};
\draw (-1,-2) node {\scriptsize $a_{-\!1}^2v_{0}, a_{-\!2}v_{0}$};
\draw (-3,-2) node {\scriptsize $a_{-\!1}v_{-\sqrt{2}}$};
\draw (1.2,-3) node {\scriptsize \dots};
\draw (-1,-3) node {\scriptsize $a_{-\!1}^3v_{0},a_{-\!2}a_{-\!1}v_{0}, a_{-\!3}v_{0}$};
\draw (-3.2,-3) node {\scriptsize \dots};
\draw (2.5,-4) node {\scriptsize $v_{2\sqrt{2}}$} edge[bend right,->,densely dotted] (2.5,-2) ;
\draw (1,-4) node {\scriptsize \dots};
\draw (-1,-4) node {\scriptsize \dots};
\draw (-3,-4) node {\scriptsize \dots};
\draw (-4.5,-4) node {\scriptsize $v_{-\!2\sqrt{2}}$};
\draw[thin, dashed] (-5.3,-5) parabola bend (-1,0) (3.3,-5); 
\draw[densely dashed,->] (-5,-4.5) -- (3,-4.5) node[anchor = north west] {$a_0$};
\draw (-4.5,-5) node {\scriptsize $-\!2\sqrt{2}$};
\draw (-3,-5) node {\scriptsize $-\!\sqrt{2}$};
\draw (-1,-5) node {\scriptsize $0$};
\draw (1,-5) node {\scriptsize $\sqrt{2}$};
\draw (2.5,-5) node {\scriptsize $2\sqrt{2}$};
\draw[densely dashed,->] (4,-5) -- (4,0) node[anchor = south] {$L_0$};
\draw (3.7,0) node {\scriptsize $0$};
\draw (3.7,-1) node {\scriptsize $1$};
\draw (3.7,-2) node {\scriptsize $2$};
\draw (3.7,-3) node {\scriptsize $3$};
\draw (3.7,-4) node {\scriptsize $4$};
\draw (8,0) node {$v_0$};
\draw (10,0) node {$v_{\sqrt{2}}$};
\draw (8,-1) node {$a_{-\!1}v_0$};
\draw (10,-1) node {$a_{-\!1}v_{\sqrt{2}}$};
\draw (6,-2) node {\scriptsize $v_{-\!\sqrt{2}}$};
\draw (7.8,-2) node {\scriptsize $a_{-\!1}^2v_0,a_{-\!2}v_0$};
\draw (10.2,-2) node {\scriptsize $a_{-\!1}^2v_{\sqrt{2}},a_{-\!2}v_{\sqrt{2}}$};
\draw (12,-2) node {\scriptsize $v_{2\sqrt{2}}$};
\draw (6,-3) node {\scriptsize $a_{-\!1}v_{-\!\sqrt{2}}$};
\draw (8,-3) node {\scriptsize \dots};
\draw (10,-3) node {\scriptsize \dots};
\draw (12,-3) node {\scriptsize $a_{-\!1}v_{2\sqrt{2}}$};
\draw (6,-4) node {\scriptsize \dots};
\draw (8,-4) node {\scriptsize \dots};
\draw (10,-4) node {\scriptsize \dots};
\draw (12,-4) node {\scriptsize \dots};
\draw[thin, dashed] (4.4,-5) parabola bend (9,0.25) (13.6,-5); 
\draw[densely dashed,->] (5,-4.5) -- (12.5,-4.5) node[anchor = north west] {$a_0$};
\draw (6,-5) node {\scriptsize $-\!\sqrt{2}$};
\draw (8,-5) node {\scriptsize $0$};
\draw (10,-5) node {\scriptsize $\sqrt{2}$};
\draw (12,-5) node {\scriptsize $2\sqrt{2}$};
\end{tikzpicture} 
 \caption{\small The basic vectors with the lowest $L_0$ grading. The left part correspond to the vacuum representation of $V_{\sqrt{2}\mathbb{Z}}$, the right part correspond to the vacuum representaion of $\calU$. Dotted curved arrows shows the shift of the $L_0$ grading to $L_0-\frac{1}{\sqrt{2}}a_0$. } \label{fig_01}
\end{center}
\end{figure}
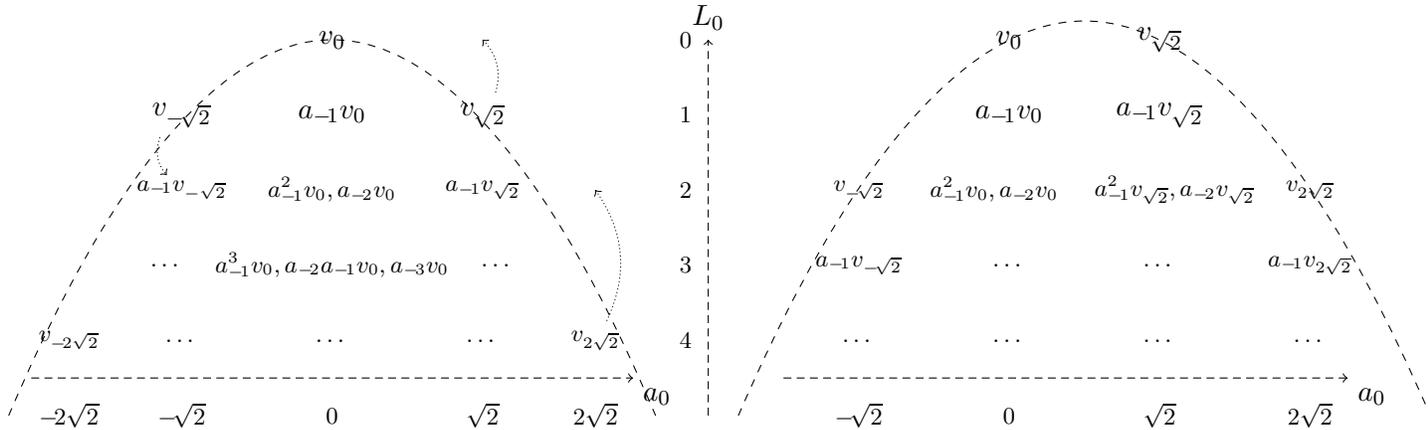

 Their characters have the form
\begin{equation}
\begin{aligned}
\chi(U_0)&=\left.\mathrm{Tr}q^{L_0}\right|_{U_0}=\sum_{k \in \mathbb{Z}} \frac{q^{k^2-k}}{(q)_\infty}=2+2q+6q^2+8q^3+\cdots,
\\ 
\chi(U_1)&=\left.\mathrm{Tr}q^{L_0}\right|_{U_1}=\sum_{k \in \mathbb{Z}+\frac12} \frac{q^{k^2-k}}{(q)_\infty}=q^{-1/4}+3q^{3/4}+4q^{7/4}+\cdots. 
\end{aligned}\label{eq:chU}
\end{equation}
\begin{Remark} \label{Re:chiU}
It is interesting to note that $\chi(U_0)=q^{-1/4}\chi(\calL_{1,1})$ and $\chi(U_1)=q^{-1/4}\chi(\calL_{0,1})$.
\end{Remark}

\subsection{Main theorem.}
Denote by $\mathrm{V}_{\Delta,c}$ the Verma module of the Virasoro algebra generated by the highest weight vector $v$
\[L_n v=0,\; \text{for}\; n>0\qquad L_0 v=\Delta v,\; Cv=cv.\]
By $\mathbb{L}_{{\Delta},c}$ denote its irreducible quotient. 
It is convenient to parametrize $\Delta$ and $c$ as 
\begin{align}
\Delta=\Delta(P,b)=\frac{(b^{-1}+b)^2}4-P^2, \qquad c=c(b)=1+6(b^{-1}+b)^2 \label{eq:Pb}
\end{align}
We denote the corresponding irreducible representation as $\mathbb{L}_{\textsc{p},b}$ keeping  in mind this parametrization. Let $b$ be generic i.e. $b^2 \not \in \mathbb{Q}$. Let $m,n \in \mathbb{Z}$ and
\begin{align} \label{eq:Pmn}
P_{m,n}=(mb^{-1}+nb)/2.
\end{align}
Then for $P\not \in \{P_{m,n}\}$ or $P=P_{m,n}, mn \leq 0$ the representation $\mathbb{L}_{\textsc{p},b}$ is isomorphic to the Verma module and has the character
$\left.\text{Tr}q^{L_0}\right|_{\mathbb{L}_{\textsc{p},b}}=q^{\Delta}/(q)_\infty$ (for the reference see e.g. \cite{Feigin Fuchs 1990} or \cite{Iohara Koga}). For $P=P_{m,n}$, $mn>0$ the Verma module contains the singular vector on the level $mn$. The irreducible representation is a quotient of Verma module by the submodule generated by this singular vector. The character of this module (which we denote by $\mathbb{L}_{(m,n)}^b$ instead of $\mathbb{L}_{\textsc{p}_{m,n},b}$) reads
\[\chi_{m,n}^b(q)=\left.\mathrm{Tr}q^{L_0}\right|_{\mathbb{L}_{(m,n)}^b}=q^{\Delta}\frac{1-q^{mn}}{(q)_\infty} \]

If $P=P_{1,1}$ then $\Delta=0$ and the corresponding irreducible  representation $\mathbb{L}_{(1,1)}^b$ has the structure of a conformal vertex algebra. The conformal vector is $L_{-2}v$ and the stress--energy tensor has the form $T(z)=\sum_n L_n z^{-n-2}$.  

The vector spaces $\mathbb{L}_{\textsc{P},b}$ are the representations of this vertex algebra. We will denote this algebra by $\mathcal{M}_b$ and denote by $T_b$ its stress--energy tensor

\begin{Theorem} \label{th:Th1}
 a) Let $b_1=b/\sqrt{1-b^2}$, $b_2=\sqrt{b^2-1}$, $b$ is generic. Then the vector space \[\mathcal{A}_b = \bigoplus_{n \in 2\mathbb{N}-1} \mathbb{L}_{(1,n)}^{b_1}\otimes \mathbb{L}_{(n,1)}^{b_2}\] has a structure of the vacuum module of a conformal vertex algebra with the stress--energy tensor $T_{b_1}+T_{b_2}$. 
 
 b) This algebra is isomorphic to $\mathcal{U}\otimes \mathcal{M}_b$, with the stress--energy tensor $T_U+T_b$.
\end{Theorem}
\begin{Remark}
The parameters $b_1$ and $b_2$ from Theorem \ref{th:Th1} satisfy the relation
\[b_1^2+b_2^{-2}=-1\] 
This relation was already mentioned in the Introduction.
\end{Remark}
\begin{Remark}
The direct sum $\bigoplus_{n\in 2\mathbb{N}} \mathbb{L}_{(1,n)}^{b_1}\otimes \mathbb{L}_{(n,1)}^{b_2}$ form a representation of the $\mathcal{A}_b$. This remark can be viewed as an addition to the part a) of the Theorem \ref{th:Th1}. One can say that the whole sum $\bigoplus_{n\in \mathbb{N}} \mathbb{L}_{(1,n)}^{b_1}\otimes \mathbb{L}_{(n,1)}^{b_2}$ is an operator algebra but nonlocal in sense that the fractional powers like $(z-w)^{1/2}$ appears in the OPE.
\end{Remark}
The proof of the Theorem \ref{th:Th1} will be given in the next section. But we can make two simple checks in advance. 

First, we can check the corollary of the Theorem \ref{th:Th1} that the central charges of the algebras $\mathcal{A}_b$ and $\mathcal{U}\otimes \mathcal{M}_b$ coincide
\[c_1+c_2=\left(1+6(\frac1{b_1}+b_1)^2\right)+ \left(1+6(\frac1{b_2}+b_2)^2\right)= -5+\left(1+6(\frac1{b}+b)^2\right)=c_{\mathcal{U}}+c\]

Second, the character of the vacuum representation of the algebra $\calA_b$ equals to
\begin{align}\label{eq:chAb}
\chi(\mathcal{A}_b)=\sum_{n\in 2\mathbb{N}-1} \chi_{1,n}^{b_1}\cdot \chi_{n,1}^{b_2}= \sum_{n\in 2\mathbb{N}-1} q^{\frac{(n-2)^2-1}{4}}\frac{(1-q^n)^2}{(q)_\infty^2}=\chi(U_0)\cdot\frac{1-q}{(q)_\infty},
\end{align}
where the last equality is an easy combinatorial statement. The equality \eqref{eq:chAb} can be rewritten as $\chi(\mathcal{A}_b)=\chi(U_0)\cdot \chi_{1,1}^b$ i.e. equality of characters of vacuum modules from the Theorem \ref{th:Th1}. Similarly
\[\chi\left(\bigoplus_{n\in 2\mathbb{N}} \mathbb{L}_{(1,n)}^{b_1}\otimes \mathbb{L}_{(n,1)}^{b_2}\right)=\chi(U_1)\cdot \chi_{1,1}^b\]

Now we consider the representations of the algebra $\calA_b$. Taking into account the isomorphism $\calA_b=\calU \otimes \calM_b$ we have tensor product representations $U_0\otimes \mathbb{L}_{\textsc{P},b}$ and $U_1\otimes \mathbb{L}_{\textsc{P},b}$. The~next theorem describes the action of $T_1(z)$ and $T_2(z)$ on these modules
\begin{Theorem}\label{th:Rep}
 Let $P\not \in \{P_{m,n}\}$, $P_1= \sqrt{b^{-1}/(b^{-1}-b)}P$, $P_2=\sqrt{b/(b-b^{-1})}P$, $i=0,1$. Then the modules $U_i\otimes \mathbb{L}_{\textsc{P},b}$ have the following decomposition with respect to the subalgebra $\calM_{b_1}\otimes \calM_{b_2}$ 
\begin{equation}
U_i\otimes \mathbb{L}_{\textsc{p},b} = \bigoplus
_{k\in \mathbb{Z}+\frac{1-i}2} \mathbb{L}_{\left(\textsc{p}_1+k b_1\right),b_1} \otimes \mathbb{L}_{\left(\textsc{p}_2+kb_2^{-1}\right),b_2}. 
\label{eq:RepAb}
\end{equation}
\end{Theorem}

One can easily check that the characters of the modules on the left side and on the right side of \eqref{eq:RepAb} are equal.


\subsection{Minimal models.} \label{Sub:Min}
In the	previous subsection we have considered the general values of the Virasoro central charge $c$. 
Now let $b_{p/p'}^2=-p/p'$, $p,p' \in \mathbb{N}$, $(p,p')=1$. The central charge equals $c_{p/p'}=1-6\dfrac{(p-p')^2}{pp'}$.

The singular values of $P$ (defined in \eqref{eq:Pmn}) possess the symmetry $P_{m,m'}=-P_{p-m,p'-m'}$. Thus we have an equality of the highest weights
\[\Delta(P_{m,m'},b_{p/p'})=\Delta(P_{p-m,p'-m'},b_{p/p'}).\] 
Therefore for $0< m < p$, $0<m'< p'$ the corresponding Verma module $\mathrm{V}_{\Delta,c}$ contains two singular vectors on the level $m m'$ and $(p-m)(p'-m')$. For example, the irreducible module $\bbL_{(1,1)}^{b}$ is the quotient of the Verma module $\rmV_{0,c}$ by the submodule generated by singular vectors of the level 1 and $(p-1)(p'-1)$. 

The module $\bbL_{(1,1)}^{b_{p/p'}}$ has the structure of the vacuum module of the conformal vertex algebra which is called the \emph{minimal model}. We will denote it by $\calM_{p/p'}$. This vertex algebra is rational, the only representations are  $\bbL_{(m,m')}^{b_{p/p'}}$ for $0<m< p$, $0<m'<p'$ with the identification $\bbL_{(m,m')}^{b_{p/p'}}=\bbL_{(p-m,p'-m')}^{b_{p/p'}}$ mentioned before. We will denote such representations as $\bbL_{(m,m')}^{p/p'}$

We want to state analogues of the Theorems \ref{th:Th1} and \ref{th:Rep} for the minimal models. Note that if $b^2=-p/p'$ then $b_1^2=-p/(p+p')$ and $b_2^2=-(p+p')/p'$ i.e. they correspond to the minimal models $\calM_{p/(p+p')}$ and $\calM_{(p+p'/p)}$.
\begin{Theorem} \label{th:MinMod}
a) The vector space 
\[ \bigoplus_{\substack{0< n <p+p' \\ n \equiv 1\bmod 2 }} \mathbb{L}_{(1,n)}^{p/(p+p')}\otimes \mathbb{L}_{(n,1)}^{(p+p')/p'}\] has a structure of the vacuum module of the conformal vertex algebra with the stress--energy tensor $T_{p/(p+p')}+T_{(p+p')/p'}$.  This algebra is isomorphic to $\mathcal{U}\otimes \mathcal{M}_{p/p'}$, with the stress--energy tensor $T_U+T_{p/p'}$.\\
b) The algebra $\mathcal{U}\otimes \mathcal{M}_{p/p'}$ has natural representations $U_i\otimes \bbL^{p/p'}_{(m,m')}$ where $i=0,1$. These modules have the following decomposition
\begin{align*}
U_i\otimes \bbL^{p/p'}_{(m,m')} \cong \!\!\!\!\!\! \bigoplus_{\substack{0< n < p+p' \\ n \equiv m+m'+i-1\bmod 2 }} \!\!\!\!\!\!\mathbb{L}_{(m,n)}^{p/(p+p')}\otimes \mathbb{L}_{(n,m')}^{(p+p')/p'}.
\end{align*}
\end{Theorem}


Denote by $\chi_{(m,m')}^{p/p'}$ the character of the irreducible module $\bbL_{(m,m')}^{p/p'}$. The following combinatorial identities follows from the \ref{th:MinMod}
\begin{align*}
\chi(U_i)\cdot \chi^{p/p'}_{(m,m')} = \!\!\!\!\!\!\!\!\!\! \sum_{\substack{0< n <p+p' \\ n \equiv m+m'+i-1\bmod 2 }} 
\!\!\!\!\!\!\!\!\!\! \chi_{(m,n)}^{p/(p+p')}\cdot \chi_{(n,m')}^{(p+p')/p'}.
\end{align*} 
The Theorem \ref{th:MinMod} has a remarkable particular case. Let $(p,p')=(2,3)$ then the central charge $c_{2/3}=0$. The minimal model $\calM_{2/3}$ has only one representation $\bbL_{(1,1)}^{2/3}$ which is trivial representation. Therefore factor $\calM_{2/3}$ can be omitted.
\begin{Theorem}\label{th:Ur}
The Urod algebra $\mathcal{U}$ has the subalgebra $\calM_{2/5}\otimes\calM_{5/3}$. The representations $U_0$  and $U_1$ have the decomposition
\[U_0=\left(\bbL_{(1,1)}^{2/5}\otimes \bbL_{(1,1)}^{5/3}\right) \bigoplus \left(\bbL_{(1,3)}^{2/5}\otimes \bbL_{(3,1)}^{5/3}\right),\qquad
U_1=\left(\bbL_{(1,2)}^{2/5}\otimes \bbL_{(2,1)}^{5/3}\right) \bigoplus \left(\bbL_{(1,4)}^{2/5}\otimes \bbL_{(4,1)}^{5/3}\right).\]
\end{Theorem}

Two commuting Virasoro algebras can be constructed explicitly in terms of Heisenberg algebra
\begin{align*}
T_{2/5}&=-\frac{1}{10\epsilon}e^{-\sqrt{2}\varphi}+\frac{1}{5} (\partial\varphi)^2+\frac{3}{5\sqrt{2}}(\partial^2\varphi)+\frac{12\epsilon}{5} (\partial\varphi)^2e^{\sqrt{2}\varphi}+\frac{3\sqrt{2}\epsilon}{5} (\partial^2\varphi)e^{\sqrt{2}\varphi}-\frac{12\epsilon^2}5e^{2\sqrt{2}\varphi},\\
T_{5/3}&=\frac{1}{10\epsilon}e^{-\sqrt{2}\varphi}+ \frac{3}{10}(\partial\varphi)^2+\frac{2}{5\sqrt{2}}(\partial^2\varphi)-\frac{2\epsilon}{5} (\partial\varphi)^2e^{\sqrt{2}\varphi} +\frac{2\sqrt{2}\epsilon}5  (\partial^2\varphi)e^{\sqrt{2}\varphi}+\frac{12\epsilon^2}5e^{2\sqrt{2}\varphi}.
\end{align*}
Direct calculation shows that $T_{2/5}$ and $T_{5/3}$ commute and satisfy \eqref{eq:T(z)OPE} with the central charges $c_{2/5}=-\frac{22}5$ and $c_{5/3}=-\frac{3}5$ correspondingly. It is clear that
\[ T_{\mathcal{U}}=T_{2/5}+T_{5/3}. \]
\begin{Remark} One can try to find the general stress--energy tensor $T(z)$ in ansatz
\[T(z)=\alpha e^{-\sqrt{2}\varphi} + \beta_1 (\partial\varphi)^2 + \beta_2(\partial^2\varphi) + \gamma_1(\partial\varphi)^2e^{\sqrt{2}\varphi} + \gamma_2(\partial^2\varphi)e^{\sqrt{2}\varphi} + \delta e^{2\sqrt{2}\varphi}.\]
The solutions of the equation \eqref{eq:T(z)OPE} are $T_{\mathcal{U}}$, $T_{2/5}$, $T_{5/3}$, standard solution $T(z)=\frac12(\partial\varphi)^2 +u(\partial^2\varphi)$ and another two deformations $T^{(1)}(z)=\alpha e^{-\sqrt{2}\varphi}+\frac12(\partial\varphi)^2$ and $T^{(2)}(z)=\frac12(\partial\varphi)^2+\frac{3}{2\sqrt{2}}(\partial^2\varphi)+\delta e^{2\sqrt{2}\varphi}$. The corresponding central charges equals $c_{\calU}=-5$, $c_{2/5}=-\frac{22}5$, $c_{5/3}=-\frac{3}5$, $c=1-12u^2$, $c^{(1)}=1$, $c^{(2)}=-\frac{25}{2}$. We can conclude that formulas for $T_{2/5}$, $T_{5/3}$ above are quite distinguished. 
\end{Remark}

\begin{figure}
\begin{center}
\begin{tikzpicture}
\draw[fill] (-0.5,0) circle (1.5pt) node [anchor=north] {\small $\mathrm{v}_1$};
\draw[fill] (-1.2,-4) circle (1.5pt) node [anchor=north] {\scriptsize $L_{-2}^{2/5}\mathrm{v}_1$};
\draw[fill] (0.2,-4) circle (1.5pt) node [anchor=north] {\scriptsize  $L_{-2}^{5/3}\mathrm{v}_1$};
\draw[fill] (-0.5,-5) node {\small  $\dots\quad \dots$};
\draw[fill] (-0.5,-6.5) node {\small  $\bbL_{(1,1)}^{2/5}\otimes \bbL_{(1,1)}^{5/3}$};
\draw[thin, loosely dashed] (-2,-7) --  (-2,-3.8) -- (-0.5,0.3) -- (1,-3.8) -- (1,-7);
\draw[fill] (3.0,0.5) circle (1.5pt) node [anchor=north] {\small $\mathrm{v}_2$};
\draw[fill] (2.3,-1.5) circle (1.5pt) node [anchor=north] {\small $L_{-1}^{2/5}\mathrm{v}_2$};
\draw[fill] (3.7,-1.5) circle (1.5pt) node [anchor=north] {\small  $L_{-1}^{5/3}\mathrm{v}_2$};
\draw[fill] (2.3,-3.5) circle (1.5pt) node [anchor=north] {\scriptsize  $L_{-2}^{2/5}\mathrm{v}_2$};
\draw[fill] (3.7,-3.5) circle (1.5pt) node [anchor=north] {\scriptsize $L_{-2}^{5/3}\mathrm{v}_2$};
\draw[fill] (3,-5) node {\small  $\dots\quad \dots$};
\draw[fill] (3,-6.5) node {\small  $\bbL_{(1,2)}^{2/5}\otimes \bbL_{(2,1)}^{5/3}$};
\draw[thin, loosely dashed] (1.7,-7) --  (1.7,-1.5) -- (3.0,0.8) -- (4.3,-1.5) -- (4.3,-7);
\draw[fill] (7.5,0) circle (1.5pt) node [anchor=north] {\small $\mathrm{v}_3$};
\draw[fill] (6.8,-2) circle (1.5pt) node [anchor=north] {\small $L_{-1}^{2/5}\mathrm{v}_3$};
\draw[fill] (8.2,-2) circle (1.5pt) node [anchor=north] {\small  $L_{-1}^{5/3}\mathrm{v}_3$};
\draw[fill] (5.4,-4) circle (1.5pt) node [anchor=north] {\scriptsize $L_{{-}2}^{2/5}\mathrm{v}_3$};
\draw[fill] (6.6,-4) circle (1.5pt) node [anchor=north] {\scriptsize $L_{{-}2}^{5/3}\mathrm{v}_3$};
\draw[fill] (8,-4) circle (1.5pt) node [anchor=north] {\scriptsize $\left(L_{{-}1}^{5/3}\right)^2\!\!\!\mathrm{v}_3$};
\draw[fill] (9.6,-4) circle (1.5pt) node [anchor=north] {\scriptsize $L_{{-}1}^{2/5}L_{{-}1}^{5/3}\mathrm{v}_3$};
\draw[fill] (7.5,-5) node {\small  $\dots\quad \dots$};
\draw[fill] (7.5,-6.5) node {\small  $\bbL_{(1,3)}^{2/5}\otimes \bbL_{(3,1)}^{5/3}$};
\draw[thin, loosely dashed] (4.8,-7) --  (4.8,-1.5) -- (7.5,0.3) -- (10.5,-1.5) -- (10.5,-7);
\draw[fill] (12.5,-1.5) circle (1.5pt) node [anchor=north] {\small $\mathrm{v}_4$};
\draw[fill] (11.8,-3.5) circle (1.5pt) node [anchor=north] {\scriptsize $L_{-1}^{2/5}\mathrm{v}_4$};
\draw[fill] (13.2,-3.5) circle (1.5pt) node [anchor=north] {\scriptsize  $L_{-1}^{5/3}\mathrm{v}_4$};
\draw[fill] (12.5,-5) node {\small  $\dots\quad \dots$};
\draw[fill] (12.5,-6.5) node {\small  $\bbL_{(1,4)}^{2/5}\otimes \bbL_{(4,1)}^{5/3}$};
\draw[thin, loosely dashed] (11.2,-7) --  (11.2,-3.5) -- (12.5,-1.2) -- (13.8,-3.5) -- (13.8,-7);
\draw[densely dashed,->] (-2.5,-5) -- (-2.5,1) node[anchor = south] {$L_0$};
\draw (13.5,0.5) [loosely dotted] --  (-2.4,0.5);
\draw (-2.4,0.5)  --  (-2.6,0.5) node [anchor=east] {\scriptsize $-1/4$};
\draw (13.5,0) [loosely dotted] --  (-2.4,0);
\draw (-2.4,0) -- (-2.6,0) node [anchor=east] {\scriptsize $0$};
\draw (13.5,-1.5) [loosely dotted] --  (-2.4,-1.5);
\draw (-2.4,-1.5) -- (-2.6,-1.5) node [anchor=east] {\scriptsize $3/4$};
\draw (13.5,-2) [loosely dotted] --  (-2.4,-2);
\draw (-2.4,-2) -- (-2.6,-2) node [anchor=east] {\scriptsize $1$};
\draw (13.5,-3.5) [loosely dotted] --  (-2.4,-3.5);
\draw (-2.4,-3.5) -- (-2.6,-3.5) node [anchor=east] {\scriptsize $7/4$};
\draw (13.5,-4) [loosely dotted] --  (-2.4,-4);
\draw (-2.4,-4) -- (-2.6,-4) node [anchor=east] {\scriptsize $2$};
\draw[densely dashed,->] (13.5,-6) -- (-2.4,-6);
\end{tikzpicture} 
\caption{\small The basic vectors in $\bigoplus_{n=1}^4 \bbL_{(1,n)}^{2/5}\otimes \bbL_{(n,1)}^{5/3}$ with $L_0$ grading $\leq 2$. The vectors $\mathrm{v}_n$, $1\leq n \leq 4$ stands for the highest weight vector of $\bbL_{(1,n)}^{2/5}\otimes \bbL_{(n,1)}^{5/3}$. By $L_n^{2/5}$ and $L_n^{5/3}$ we denote the components of $T_{2/5}(z)$ and $T_{5/3}(z)$ correspondingly. \label{fig_U}}
\end{center}
\end{figure}
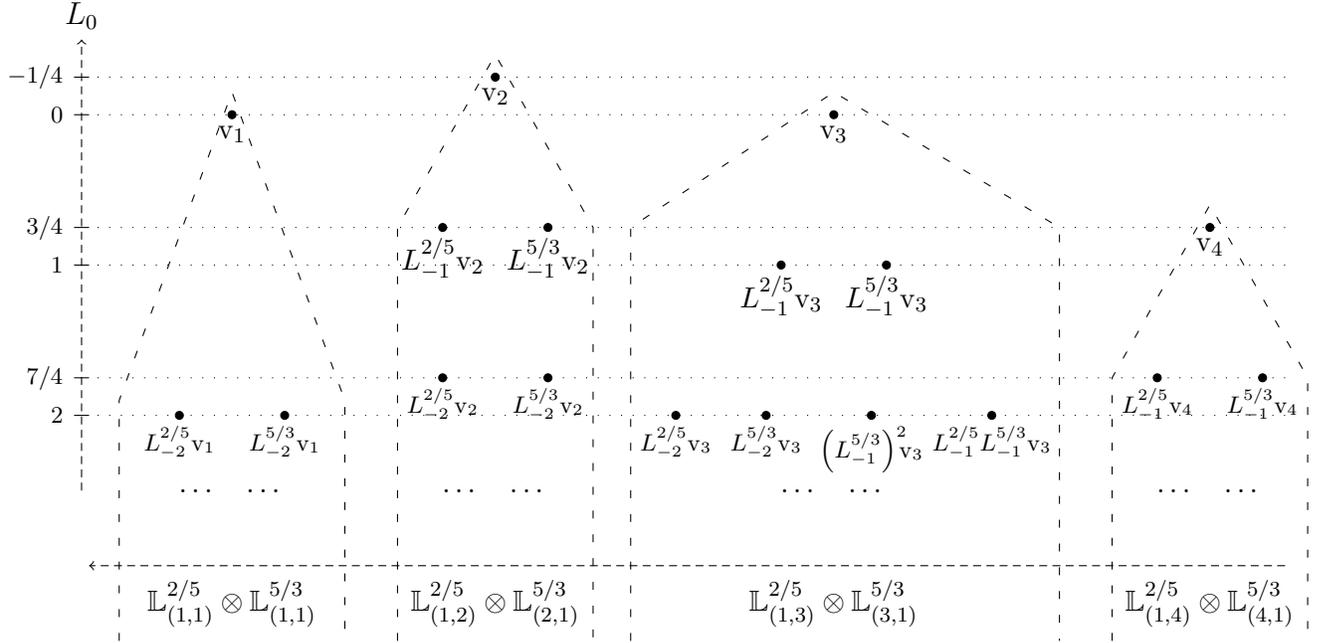

At the end of the subsection we remark two identities
\begin{equation}
\begin{aligned}
\chi(\calL_{0,1})=q^{-1/4}\left(\chi_{(1,2)}^{2/5}\cdot \chi_{(2,1)}^{5/3}+\chi_{(1,4)}^{2/5}\cdot \chi_{(4,1)}^{5/3}\right),\\ \chi(\calL_{1,1})=q^{-1/4}\left(\chi_{(1,1)}^{2/5}\cdot \chi_{(1,1)}^{5/3}+\chi_{(1,3)}^{2/5}\cdot \chi_{(3,1)}^{5/3}\right). 
\end{aligned}\label{eq:chi_sl2Vir2535}
\end{equation}
These identities link characters of the Level 1 representations $\Asl(2)$ and characters of minimal models and follow from the Theorem \ref{th:Ur} and Remark \ref{Re:chiU}. These identities will be discussed in Section \ref{Seq:Comb} from the combinatorial point of view.

\subsection{$c=-5$ description.} We can also study $\mathcal{U}$ as  the Virasoro algebra representation with the central charge $c_{\mathcal{U}}=-5$. The corresponding parameter $b_{\calU}=i\dfrac{1+\sqrt{5}}{2}$ (see parameterization \eqref{eq:Pb}) is generic by means $b_\calU^2 \not \in \mathbb{Q}$. The unique singular vector in Verma module $\rmV_{(m,n)}^{b_\calU}$ has $L_0$ grading $\Delta(P_{m,n},b_\calU)+mn=\Delta(P_{m,-n},b_\calU)$. These fact can be also written in a short exact sequence
\[0 \rightarrow \bbL_{(m,-n)}^{b_\calU} \rightarrow \rmV_{(m,n)}^{b_\calU} \rightarrow \bbL_{(m,n)}^{b_\calU} \rightarrow 0\]
We also need the projective modules  in the BGG category $\mathcal{O}$ of Virasoro representation with central charge $-5$. By $\mathcal{P}_{m,n}$ we denote the unique projective module such that $\mathrm{Hom}(\mathcal{P}_{m,n},\bbL_{m,n})\neq 0$. These projective modules have the description similar to the projective modules for the category $\mathcal{O}$ for $\mathfrak{sl}(2)$. For $mn\geq 0$ we have an isomorphism $\mathcal{P}_{(m,n)}^{b_\calU} \cong \mathbb{V}_{(m,n)}^{b_\calU}$. If $mn<0$ then $\mathcal{P}_{(m,-n)}^{b_\calU}$ is defined by a short exact sequences
\[0 \rightarrow \rmV_{(m,n)}^{b_\calU} \rightarrow \mathcal{P}_{(m,-n)}^{b_\calU} \rightarrow \rmV_{(m,-n)}^{b_\calU} \rightarrow 0.\]
The last short exact sequence also follows from the duality theorem for the category $\mathcal{O}$, see \cite{BGG} or  \cite[Theorem 1.2]{Iohara Koga}.

\begin{``Theorem''} \label{th:c=-5}
The modules $\bbL_{(1,n)}^{2/5}\otimes \bbL_{(n,1)}^{5/3}$, $1\leq n \leq 4$ have the following decomposition with respect to diagonal Virasoro $L_n=L_n^{2/5}+L_n^{5/3}$
\begin{align*}
&\bbL_{(1,1)}^{2/5}\otimes \bbL_{(1,1)}^{5/3}\cong \bigoplus_{n \in 2\mathbb{N}-1} \bbL_{(n,n)}^{b_\calU}, &&\bbL_{(1,3)}^{2/5}\otimes \bbL_{(3,1)}^{5/3}\cong \bigoplus_{n \in 2\mathbb{N}-1} \mathcal{P}_{(n,-n)}^{b_\calU},\\
&\bbL_{(1,2)}^{2/5}\otimes \bbL_{(2,1)}^{5/3}\cong \mathcal{P}_{(0,0)}^{b_\calU}\oplus \bigoplus_{n \in 2\mathbb{N}} \mathcal{P}_{(n,-n)}^{b_\calU},  &&\bbL_{(1,4)}^{2/5}\otimes \bbL_{(4,1)}^{5/3}\cong \bigoplus_{n \in 2\mathbb{N}} \mathbb{L}_{(n,n)}^{b_\calU}.
\end{align*}
\end{``Theorem''}
Due to Theorem \ref{th:Ur} the first row gives the decomposition of $U_0$ and the second row gives the decomposition of $U_1$.

We will not prove this fact in this paper (and therefore we called it ``Theorem'').

\section{Constructions and Proofs} \label{Se:proof}

\subsection{Quantum Hamiltonian reduction of $\Asl(2)$.} For the reader convenience we first recall the quantum Hamiltonian (or the Drinfeld--Sokolov) reduction of $\Asl(2)$, see e.g. \cite{FrenkelBook} for the reference. We denote by $\calV_{h,k}$ the Verma module of the $\Asl(2)$, and by $\calL_{h,k}$ its irreducible quotient.

Let $V$ be a representation of vertex algebra $\calL_{0,k}$. It is convenient to consider $\Asl(2)$ generators as modes of the fields $e(z), f(z), h(z)$
\begin{align*}
e(z)=\sum_{n\in \mathbb{Z}} e_n z^{-n-1}, \quad h(z)=\sum_{n\in \mathbb{Z}} h_n z^{-n-1}, \quad f(z)=\sum_{n\in \mathbb{Z}} f_n z^{-n-1} 
\end{align*}
The relations of the $\Asl(2)_k$ can be written in terms of OPE\begin{align*}
h(z)e(w) \sim \frac{2 e(w)}{(z-w)},\quad h(z)f(w) \sim -\frac{2f(w)}{(z-w)}, \quad h(z)h(w)\sim \frac{2 k}{(z-w)^2} \\
e(z)f(w)\sim \frac{k}{(z-w)^2}+\frac{h(w)}{(z-w)},\quad e(z)e(w)\sim 0, \quad f(z)f(w) \sim 0
\end{align*}
The stress--energy tensor is given by the Sugavara formula \[T_{\mathrm{Sug}}(z)=\frac1{2(k+2)}:\left(\frac12 h^2(z)+e(z)f(z)+f(z)e(z)  \right):\]
The central charge of $T_{\mathrm{Sug}}$ equals $c_{\mathrm{Sug}}=\frac{3k}{k+2}$. 

Introduce the anticommuting operators $\psi_n, \psi^*_n$ with the relations \[\psi(z)=\sum_{n}\psi_n z^{-n}, \psi^*(z)=\sum_{n}\psi^*_n z^{-n},\qquad \psi(z) \psi^*(w)\sim \frac1{z-w},\quad \psi(z) \psi(w)\sim\psi^*(z) \psi^*(w)\sim0\]
Operators $\psi_n, \psi^*_n$ generates Clifford algebra. By $\Lambda$ we denote the Fock representation generated by the vector $v$
\[\psi_nv=0\;\, \text{for $n \geq 0$},\quad  \psi^*_nv=0\;\, \text{for $n>0$}.\]
We introduce the grading on $\Lambda$ by $\deg(v)=0, \deg(\psi_n)=1, \deg(\psi^*_n)=-1$. The operator
\[\mathrm{Q} = \oint_{|z|=1}(e(z)+1)\psi(z)dz \]
acts on the space $V \otimes \Lambda$. It is easy to see that $\mathrm{Q}^2=0$. We denote by $\mathrm{H}^i_{\mathrm{DS}}(V)$ the cohomology of the complex $(V\otimes \Lambda,\mathrm{Q})$, where $i$ stands for the grading on $\Lambda$. These cohomology are called quantum Hamiltonian (or the Drinfeld--Sokolov) reduction of the~$V$.

The following proposition is standard (see e.g. \cite[Sec. 15.1.8]{FrenkelBook})

\begin{Proposition}\label{pr:DSvac}
If $V$ has a structure of vertex algebra then $\mathrm{H}_{\mathrm{DS}}(V\otimes \Lambda,\rmQ)$ has a structure of~vertex algebra.
\end{Proposition}

This proposition basically follows from the formula $[\rmQ,Y(v;z)]=Y(\rmQ v;z)$ (see \cite[Corollary 3.3.8]{FrenkelBook}). In other words action of $\rmQ$ on the space of states $v$ is equivalent to the commutator with $\rmQ$ on the space of corresponding operators $Y(v;z)$.

\begin{Theorem}[e.g. {\cite[Ch. 15]{FrenkelBook}}] Let $V \cong \calL_{0,k}$, $k \neq -2$.  Then $\mathrm{H}^i_{\mathrm{DS}}(\calL_{0,k})=0$ for $i \neq 0$ and $\mathrm{H}^0_{\mathrm{DS}}(\calL_{0,k})$ is isomorphic to the vacuum representation of the algebra Virasoro with the central charge $c$, where $c=c(b_{{k{+}2}/1})$,  $b_{{k{+}2}/1}=\sqrt{-(k+2)}$. 
\end{Theorem}

The stress--energy tensor which generates the Virasoro symmetry reads
\[ T_{\mathrm{DS}}(z)=T_{\mathrm{Sug}}(z)+\frac12 \partial_z h(z)-\psi(z)\partial\psi^*(z).
\]
The operator $T_{\mathrm{DS}}(z)$ commutes with $\rmQ$ (and corresponds to the vector in the cohomology $\mathrm{H}^0_{\mathrm{DS}}(\calL_{0,k})$). Therefore, this Virasoro algebra acts on $\mathrm{H}_{\mathrm{DS}}(V)$ if $V$ is an  any level $k$ representation of $\Asl(2)$.

 We say that the pair $(h,k)$ is \emph{generic} if there is no $m,n \in \mathbb{Z}_{\geq 0}$ such that $h+m(k+2)=n$ or $k-h+m(k+2)=n$. These conditions are equivalent to the fact that Shapovalov form on the Verma module $\calV_{h,k}$ is non degenerate \cite{KacKazhdan}. Therefore the pair $(h,k)$ is generic if and only if the Verma module $\calV_{h,k}$ is irreducible. 
 
Another important example is generic $k$ and integer $h=n\in\mathbb{Z}_{\geq 0}$. In this case the irreducible module $\calL_{h,k}$ is integrable with respect to $\mathfrak{sl}(2)$ generated by $e,h,f$.

\begin{Theorem}[\cite{Feigin:1990}] \label{th:DS} Let $b=b_{(k{+}2)/1}=\sqrt{-(k+2)}$.

a) Let $n \in \mathbb{Z}_{\geq 0}$, $V \cong \calL_{n,k}$. Then $\mathrm{H}^i_{\mathrm{DS}}(\mathcal{L}_{h,k})=0$, for $i \neq 0$ and $\mathrm{H}^0_{\mathrm{DS}}(\calL_{n,k})=\bbL_{(n{+}1,1)}^{b}$.

b) Let $(h,k)$ be generic, $V \cong \calL_{h,k}$. Then $\mathrm{H}^i_{\mathrm{DS}}(\mathcal{L}_{h,k})=0$, for $i \neq 0$ and $\mathrm{H}^0_{\mathrm{DS}}(\calL_{h,k})=\bbL_{\textsc{p},b}$, where $P=\frac{b^{-1}(h+1)+b}2$.
\end{Theorem}

\subsection{Coset.}  Assume that $k \not \in \mathbb{Q}$.

Consider the tensor product of two $\Asl(2)$ modules $\calL_{i,1}\otimes \calL_{h,k}$, where $i=0,1$. There is an action of the algebra $\Asl(2)\otimes \Asl(2)$ on this space, we denote by $e^{(1)}_n,h^{(1)}_n,f^{(1)}_n$ the generators of the first factor and by $e^{(2)}_n,h^{(2)}_n,f^{(2)}_n$ the generators of the second factor. 

The space $\calL_{i,1}\otimes \calL_{h,k}$ becomes a level $k+1$ representation under the diagonal action of $\Asl(2)^\Delta$: $e^{\Delta}_n=e^{(1)}_n+e^{(2)}_n, h^{\Delta}_n=h^{(1)}_n+h^{(2)}_n, f^{\Delta}_n=f^{(1)}_n+f^{(2)}_n$. It was noticed in \cite{GKO} that there is a Virasoro algebra which commute with an action of $\Asl(2)^\Delta$
\[T_{\mathrm{Coset}}=T_{\mathrm{Sug}}^{(1)}+T_{\mathrm{Sug}}^{(2)}-T_{\mathrm{Sug}}^{\Delta} \]
This Virasoro algebra is called coset Virasoro algebra, its central charge equals $c=c\left(b_{(k+2)/(k+3)}\right)$, where $b_{(k+2)/(k+3)}=\sqrt{-\frac{k+2}{k+3}}$.

\begin{Theorem}\label{th:coset} Let $b=b_{(k+2)/(k+3)}$.

a) The tensor product of the vacuum modules have the decomposition as a $\Vir\oplus \Asl(2)^\Delta$ module
\begin{align*}
\calL_{0,1}\otimes \calL_{0,k} = \bigoplus_{n \in \mathbb{Z}_{\geq 0}} \bbL_{(1,2n+1)}^{b_{(k{+}2)/(k{+}3)}} \otimes \calL_{2n,k{+}1},
\end{align*}

b) Let $(h,k)$ be generic, $P=\frac12\frac{(1+h)}{\sqrt{-(k+2)(k+3)}}$. Then we have a decomposition of the $\calL_{i,1}\otimes \calL_{h,k}$ as a $\Vir\oplus \Asl(2)^\Delta$ module
\begin{align*}
\calL_{i,1}\otimes \calL_{h,k} = \bigoplus_{n \in \mathbb{Z}+\frac{i}2} \bbL_{\textsc{p}+nb,b} \otimes \calL_{h+2n,k{+}1}, 
\end{align*}
\end{Theorem}
This theorem seems to be standard but the authors
could not find a good reference. We give a sketch of the proof. 
\begin{proof}[Sketch of the proof] The proof is based on the same two arguments as the proof in the more difficult case of admissible representations, see Theorem \ref{th:DSmm} (and see \cite[Theorem 10.2]{Iohara Koga} as the reference). Namely first one checks the identity of characters of the module on the left side and right side. Second one checks that there are no extensions between modules on the right side.
\end{proof}

Now we can prove Theorem \ref{th:Th1} a). We construct vertex algebra algebra $\calA_b$ by the  Drinfeld--Sokolov reductions to the space $V=\calL_{0,1}\otimes \calL_{0,k}$, where $b=\sqrt{-(k+2)}$. Using the Theorems \ref{th:DS} and \ref{th:coset} we get
\begin{align}
\mathrm{H}_{\mathrm{DS}}^\Delta\left(\calL_{0,1}\otimes \calL_{0,k}\right)=\bigoplus_{n \in \mathbb{Z}_{\geq 0}}  \bbL_{(1,2n+1)}^{b_{k{+}2/k{+}3}} \otimes  \mathrm{H}_{\mathrm{DS}} (\calL_{2n,k{+}1})=\bigoplus_{n \in \mathbb{Z}_{\geq 0}}  \bbL_{(1,2n{+}1)}^{b_{k{+}2/k{+}3}} \otimes  \bbL_{(2n{+}1),1}^{b_{{k{+}3}/1}} \label{eq:DSprod}
\end{align}
This space has a vertex algebra structure due to the Proposition \ref{pr:DSvac}. If we put $b=\sqrt{-(k+2)}$, then $b_1=b_{k{+}2/k{+}3}$,  $b_2=b_{k{+}3/1}$ and we obtain the statement of the Theorem \ref{th:Th1} a).

The total stress--energy tensor is given by the formula
\begin{multline*}
T_{\mathrm{tot}}(z)=T_{\mathrm{Coset}}(z)+T_{\mathrm{DS}}^{\Delta}(z)=T_{\mathrm{Sug}}^{(1)}(z)+T_{\mathrm{Sug}}^{(2)}(z)+\frac12 \partial_z \left(h^{1}(z)+h^{2}(z)\right)-\psi(z)\partial\psi^*(z)
\end{multline*}

\subsection{Proofs of Theorems \ref{th:Th1}b), \ref{th:Rep},  \ref{th:MinMod}.}
\begin{proof}[Proof of Theorem \ref{th:Th1} b).]
 Recall that level 1 representations $\calL_{i,1}$ can be constructed in terms of one  field $\varphi(z)$ (see subsection \ref{subU}).

Introduce the differential $\mathrm{Q}_\epsilon$
\[\mathrm{Q}_\epsilon= \oint_{|z|=1}(\epsilon e^{(1)}(z)+e^{(2)}(z)+1)\psi(z) \]
If $\epsilon=1$ then $\rmQ_1$ coincides with the Drinfeld--Sokolov differential for $\Asl(2)^\Delta$. Actually for any $\epsilon \neq 0$  we can rescale $e^{(1)}$ and $f^{(1)}$ by $\epsilon$ and again get the Drinfeld--Sokolov differential for $\Asl(2)^\Delta$. But for $\epsilon=0$ the situation is changed and $\rmQ_0$ becomes the  Drinfeld--Sokolov differential for $\Asl(2)^{(2)}$.

Therefore the cohomology space of $\rmQ_0$ has the form
$\mathrm{H}_{\mathrm{DS}}^{(2)}\left(\calL_{0,1}\otimes \calL_{0,k}\right)=\calL_{0,1}\otimes  \bbL_{(1,1)}^b$, where $b=\sqrt{-(k+2)}$. If we decompose $T_{\mathrm{tot}}=T_{\calL_{0,1}}+T_{\mathrm{DS}}^{(2)}$ then we get 
\begin{equation}\label{eq:L01}
T_{\calL_{0,1}}=T_{\mathrm{tot}}-T_{\mathrm{DS}}^{	(2)}=\frac12(\partial\varphi)^2+\frac{1}{\sqrt{2}}\partial^2\varphi.
\end{equation}

We can compute the cohomology of $\rmQ_\epsilon=\rmQ_0+\epsilon\oint e^{(1)} \psi(z)$ by use of spectral sequence. The $E^1$ term of this sequence equals to the cohomology of $\rmQ_0$. Since the Drinfeld--Sokolov cohomology $\mathrm{H}_{\mathrm{DS}}^i$ vanishes for $i \neq 0$ then the spectral sequence degenerates in $E^1$ term. 

Thus $\calL_{0,1}\otimes  \bbL_{(1,1)}^b$ is isomorphic to $\bigoplus_{n \in \mathbb{Z}_{\geq 0}}  \bbL_{(1,2n{+}1)}^{b_{k{+}2/k{+}3}} \otimes  \bbL_{(2n{+}1),1}^{b_{{k{+}3}/1}}$ as a $L_0$ graded vector space. But the representative of the cohomology classes are deformed like $A \rightarrow A+\epsilon A_1+e^2A_2+\dots$. The Theorem \ref{th:Th1} b) states that $E^\infty=\bigoplus_{n \in \mathbb{Z}_{\geq 0}}  \bbL_{(1,2n{+}1)}^{b_{k{+}2/k{+}3}} \otimes  \bbL_{(2n{+}1),1}^{b_{{k{+}3}/1}}$ is isomorphic to $E^1=\calL_{0,1}\otimes  \bbL_{(1,1)}^b$ as the vertex  algebra but has a different chiral structure. 

We will prove the isomorphism of these vertex algebras by direct calculation \footnote{Probably it would be better to deduce this fact from the vanishing of the space  infinitesimal deformations of $\calL_{0,1}\otimes  \bbL_{(1,1)}^b$}.  

In order to study the deformation of the representative of cohomology classes we denote $\partial \chi = \psi \psi^*-\frac12 h^{(2)}.$ The field $\partial \chi$ has the properties 
\[[\mathrm{Q}_0,\partial\chi(z)]=-\psi(z) \qquad \partial\chi(z)\partial\chi(w)=\frac{k+2}{2(z-w)^2} \]
We give formulas for the representatives of the $\rmQ_\epsilon$ cohomology space. Denote
$\partial\tilde{\varphi}=\partial\varphi-\epsilon\sqrt{2} \partial\chi e^{\sqrt{2}\varphi}$. Then 
\[
[\mathrm{Q}_\epsilon,\partial\tilde{\varphi}(z)]=0, \quad \partial\tilde{\varphi}(z)\partial\tilde{\varphi}(w)=\frac{1}{(z-w)^2}.
\]
In other words the $\partial\tilde{\varphi}$ is free bosonic field in the cohomology of $\mathrm{Q}_\epsilon$. The stress--energy tensor of the $\bbL_{(1,1)}^b$ which we denote by $T_{\mathrm{DS}}^{(2)}$ is deformed to
\[
\tilde{T}=T_{\mathrm{DS}}^{(2)}+\epsilon \left( \sqrt{2}\cdot\partial\chi\partial\varphi -\frac{k{+}1}{2}\cdot \left(2(\partial\varphi)^2 +\sqrt{2}(\partial^2\varphi)\right)\right)  e^{\sqrt{2}\varphi}-\epsilon^2\frac{k{+}2}{2}\cdot e^{2\sqrt{2}\varphi}\]
This current $\tilde{T}(z)$ commutes with $\partial\tilde{\varphi}(w)$ and satisfies $[\rmQ_\epsilon,\tilde{T}(z)]=0$ and  Virasoro OPE with the central charge $c=c(b)$. The stress--energy tensor for the field $\tilde{\varphi}$ can be found as $T_{\mathrm{tot}}-\tilde{T}$
\begin{multline*}
T_{\mathcal{U}}=\frac12(\partial\varphi)^2 +\frac1{\sqrt{2}}(\partial^2\varphi)-\epsilon \left( \sqrt{2}\cdot\partial\chi\partial\varphi -\frac{k{+}1}{2}\cdot \left(2(\partial\varphi)^2 +\sqrt{2}(\partial^2\varphi)\right)\right)  e^{\sqrt{2}\varphi}+\epsilon^2\frac{k{+}2}{2}\cdot e^{2\sqrt{2}\varphi}
\\ = \frac12(\partial\tilde{\varphi})^2 +\frac1{\sqrt{2}}(\partial^2\tilde{\varphi})+\epsilon\frac{k{+}1}{2}\cdot \left(2(\partial\tilde{\varphi})^2 +\sqrt{2}(\partial^2\tilde{\varphi})\right)  e^{\sqrt{2}\varphi} .
\end{multline*}
This last formula for $T_{\mathcal{U}}$ coincides with equation \eqref{eq:Tu} after rescaling of $\varepsilon$. 
\end{proof}

 The representations of the algebra $\calA_b$ can be constructed by use of Drinfeld--Sokolov reduction too.
 
\begin{proof}[Proof of Theorem \ref{th:Rep}]
Consider the representations of the vertex algebra $\calL_{0,1}\otimes \calL_{0,k}$ namely 
$\calL_{i,1}\otimes \calL_{h,k}$, where $i=0,1$, and $(h,k)$ is generic.  After quantum Hamiltonian reduction with respect to diagonal we get the representations of algebra $\calA_b$. On the other hand one can use Theorems \ref{th:DS} and \ref{th:coset} and obtain decomposition with respect to subalgeba $\calM_{b_1}\otimes \calM_{b_2}$
\[
\mathrm{H}_{\mathrm{DS}}^\Delta\left(\calL_{i,1}\otimes \calL_{h,k}\right)=\bigoplus_{n \in \mathbb{Z}+\frac{i}{2}} \bbL_{\tilde{\textsc{p}}_1+nb_1,b_1} \otimes \mathrm{H}_{\mathrm{DS}}^\Delta\left(\calL_{h+2n,k{+}1}\right)=\bigoplus_{n \in \mathbb{Z}+\frac{i}{2}} \bbL_{\tilde{\textsc{p}}_1+nb_1,b_1} \otimes \bbL_{\tilde{\textsc{p}}_2+nb_2^{-1},b_2}
\]
where we used notations $b=\sqrt{-(k+2)}$, $b_1=b_{(k{+}2)/(k{+}3)}$,  $b_2=b_{(k{+}3)/1}$, $\tilde{P}_1=\frac12\frac{(1+h)}{\sqrt{-(k+2)(k+3)}}$, $\tilde{P}_2=\frac{h-2-k}{\sqrt{-k-3}}$. Now we slightly change these notations \[P_1=\tilde{P}_1+b_1/2=\frac{h-k-1}{2\sqrt{-(k+3)(k+2)}},\qquad P_2=\tilde{P}_2+b_2^{-1}/2=\frac{h-k-1}{2\sqrt{-(k+3)}}. \] Then we have the following decompositions for the representations of $\calA_b$ with respect to subalgeba $\calM_{b_1}\otimes \calM_{b_2}$
\[\mathrm{H}_{\mathrm{DS}}^\Delta\left(\calL_{i,1}\otimes \calL_{h,k}\right)= \bigoplus
_{n \in \mathbb{Z}+\frac{1-i}{2}} \mathbb{L}_{\left(\textsc{p}_1+n b_1\right),b_1} \otimes \mathbb{L}_{\left(\textsc{p}_2+n b_2^{-1}\right),b_2},\]

We want to prove that these $\mathcal{A}_b$ representations are isomorphic to representations $U_i\otimes \chi(\bbL_{P,b}$, where $P=\frac{h-k-1}{\sqrt{-k-2}}$. Since the last representations are irresucible it is enough to show the equality of characters.

We use the same tool as in the proof below, namely we  compute cohomology of the operator $\mathrm{Q}_0$ (using Theorem \ref{th:DS})
$\mathrm{H}_{\mathrm{DS}}^{(2)}\left(\calL_{i,1}\otimes \calL_{h,k}\right)=\calL_{i,1}\otimes  \bbL_{P,b}$. The character of  $E^1$ term is equal to the character of $E^\infty$ term since the spectral sequence degenerates. Using the formula  
\eqref{eq:L01} we get 
\[\chi(U_i)\chi(\bbL_{P,b})=\sum
_{n \in \mathbb{Z}+\frac{1-i}{2}} \chi(\mathbb{L}_{\left(\textsc{p}_1+n b_1\right),b_1}) \chi(\mathbb{L}_{\left(\textsc{p}_2+n b_2^{-1}\right),b_2}).\]
\end{proof}

\begin{Remark} The isomorphism of conformal vertex algebras in Theorem \ref{th:Th1} can be made more explicit. This theorem states in particular that one can find in $\calU\otimes \calM_b$ two commuting Virasoro algebras with central charges $c(b_1)$ and $c(b_2)$. In other words one can express $T_{b_1}$ and $T_{b_2}$ in terms of free field $\varphi$, exponents $e^{n\sqrt{2}\varphi}$ and stress--energy tensor $T_b$. Namely
\begin{multline}
T_{b_1}=\frac{b+b^{-1}}{2(b-b^{-1})\epsilon}e^{-\sqrt{2}\varphi}+\frac{b}{2(b-b^{-1})}(\partial \varphi)^2 -\frac{b^{-1}}{\sqrt{2}(b-b^{-1})}\partial^2\varphi -\frac{(1+2b^{-2})\epsilon}{b^2-b^{-2}}(\partial\varphi)^2e^{\sqrt{2}\varphi} \\ -\frac{\sqrt{2}b^{-1}\epsilon}{b-b^{-1}}(\partial^2\varphi) e^{\sqrt{2}\varphi}  -\frac{2 \epsilon^2}{b^2-b^{-2}}e^{2\sqrt{2}\varphi} -\frac{b^{-1}}{b-b^{-1}}T_{b} -\frac{2\epsilon}{b^2-b^{-2}} T_{b}e^{\sqrt{2}\varphi}, \label{eq:Tb1}
\end{multline}
\begin{multline}
T_{b_2}=-\frac{b+b^{-1}}{2(b-b^{-1})\epsilon}e^{-\sqrt{2}\varphi} -\frac{b^{-1}}{2(b-b^{-1})}(\partial \varphi)^2 +\frac{b}{\sqrt{2}(b-b^{-1})}\partial^2\varphi +\frac{(2b^2+1)\epsilon}{b^2-b^{-2}}(\partial\varphi)^2e^{\sqrt{2}\varphi} \\ +\frac{\sqrt{2}b\epsilon}{b-b^{-1}}(\partial^2\varphi) e^{\sqrt{2}\varphi}  +\frac{2\epsilon^2}{b^2-b^{-2}}e^{2\sqrt{2}\varphi} +\frac{b}{b-b^{-1}}T_{b} +\frac{2\epsilon}{b^2-b^{-2}} T_{b}e^{\sqrt{2}\varphi}. 
\label{eq:Tb2}
\end{multline}
It is easy to see that $T_{b_1}+T_{b_2}=T_{b}+T_{\calU}$. 

If $b=b_{2/3}$ then $T_b=0$ and these formulas reduces to formulas for $T_{2/5}$ and $T_{5/3}$ from Subsection \ref{Sub:Min}. We will use formulas \eqref{eq:Tb1},\eqref{eq:Tb2} in the Section \ref{Se:Equations}.
\end{Remark}


\begin{proof}[Proof of Theorem \ref{th:MinMod}]
Proof goes similar to the proof for generic $k$ but we need to use more delicate results about the coset construction and the quantum Hamiltonian reduction. 

Recall $k\in \mathbb{C}$ is called \emph{admissible level} for $\Asl(2)$ if $k=-2+\frac{p}{p'}$, where $p,p' \in \mathbb{N}$ are coprime and $p\geq 2$. The  module $\calL_{h,k}$ is called \emph{admissible module} if $k$ is an admissible level and $h=h(m,m';k)$ for $0<m<p$, $0<m'<p'$, where $h(m,m';k)=(m-1)-(k+2)(m'-1)$. 

The following two theorems are analogues of Theorems \ref{th:DS} and \ref{th:coset} for admissible representations. 

\begin{Theorem}[{\cite[Sec. 4.]{Feigin1:1990}}]
\label{th:DSmm} Let $\calL_{h(m,m';k),k}$ be an admissible module.  Then $\mathrm{H}^i_{\mathrm{DS}}(\mathcal{L}_{h(m,m';k),k})=0$, for $i \neq 0$ and $\mathrm{H}^0_{\mathrm{DS}}(\calL_{h(m,m';k),k})=\bbL_{m,m'}^{p/p'}$. 
\end{Theorem}

\begin{Theorem}[{\cite[Theorem 10.2]{Iohara Koga}}]
\label{th:cosetmm} Let $\calL_{h(m,m';k),k}$ be an admissible module.  Then we have a decomposition of the $\calL_{i,1}\otimes \calL_{h(m,m';k);k}$ as a $\Vir\oplus \Asl(2)^\Delta$ module
\[
\calL_{i,1}\otimes \calL_{h(m,m';k),k} = \bigoplus_{\substack{0< n < p+p' \\ n \equiv m+m'+i-1\bmod 2 }} \!\!\!\bbL_{m,n}^{p/(p+p')} \otimes \calL_{h(n,m';k+1),k{+}1}.
\]
\end{Theorem}
Now we take the tensor product $\calL_{i,1}\otimes \calL_{h,k}$ for admissible $\calL_{h,k}$ and apply quantum Hamiltonian reduction with respect to diagonal $\Asl(2)$
\[
\mathrm{H}_{\mathrm{DS}}^\Delta\left(\calL_{i,1}\otimes \calL_{h(m,m';k),k}\right)=\bigoplus_{\substack{0< n < p+p' \\ n \equiv m+m'+i-1\bmod 2 }} \!\!\!\bbL_{m,n}^{p/(p+p')} \otimes \bbL_{n,m'}^{(p+p')/p'},
\]
as needed in Theorem \ref{th:MinMod}. The remaining arguments are the same as in the proof of Theorems \ref{th:Th1} and \ref{th:Rep}.
\end{proof}

\subsection{Proof of Theorem \ref{th:Ur}}
The Theorem \ref{th:Ur} follows from the more general Theorem~\ref{th:MinMod}. In this subsection we give another proof of this fact using the asymptotic dimensions.

Let $V$ be the highest weight representation of the Virasoro algebra. Recall that the triple $(A,B,C)$ is called asymptotic dimension of $V$ if
\[\mathrm{Tr}e^{-2\pi t L_0} \sim A \cdot t^B \cdot\exp\left(\frac{\pi C}{12 t}\right)\quad \text{where $t\rightarrow 0$}. \]
The $C$ is also called the effective central charge. 

It is known \cite{KacWakimoto} that the level 1 representations of $\Asl(2)$ have the effective central charge $C=1$. Therefore the modules $U_0$ and $U_1$ have the effective central charge $C=1$.

These spaces are representations of the algebra $\Vir\oplus\Vir$ due to formulas for $T_{2/5}$ and $T_{5/3}$ in Subsection \ref{Sub:Min}. The only Virasoro representations with the effective central charge less then 1 are the minimal model representations \cite{KacWakimoto}. Therefore $U_0$ and $U_1$ decompose into direct sum of the tensor products of the $(2/5)$ and $(5/3)$ minimal model representations (since there is no $\mathrm{Ext}$'s between minimal models representations).

The conformal dimensions of these minimal model representations are \begin{align*}
\Delta(P_{1,1},b_{2/5})=\Delta(P_{1,4},b_{2/5})=0, \qquad\Delta(P_{1,2},b_{2/5})=\Delta(P_{1,3},b_{2/5})=-\frac15,\\
\Delta(P_{1,1},b_{5/3})=0, \quad \Delta(P_{2,1},b_{5/3})=-\frac1{20}, \quad\Delta(P_{3,1},b_{5/3})=\frac15,\quad \Delta(P_{4,1},b_{5/3})=\frac34.
\end{align*}
The eigenvalues of $L_0$ on $U_0$ belong to $\mathbb{Z}$ and the eigenvalues on $U_1$ belong to $\mathbb{Z}-\frac14$. Therefore $U_0$ is decomposed into the sum of the representations $\bbL_{(1,1)}^{2/5}\otimes \bbL_{(1,1)}^{5/3}$, $\bbL_{(1,3)}^{2/5}\otimes \bbL_{(3,1)}^{5/3}$ and $U_0$ is decomposed into the sum of the representations $\bbL_{(1,2)}^{2/5}\otimes \bbL_{(2,1)}^{5/3}$, $\bbL_{(1,4)}^{2/5}\otimes \bbL_{(4,1)}^{5/3}$. Comparing the first terms in the $q$-expansion of the characters we get the Theorem \ref{th:Ur}.

\begin{Remark} Denote by $\Phi_{1,n} \Phi_{n,1}$ the operator which corresponds to the highest weight vector of the representation $\bbL_{(1,n)}^{2/5}\otimes \bbL_{(n,1)}^{5/3}$. In terms of the free field $\varphi$ these operators have the form
\begin{align*}
\Phi_{1,1} \Phi_{1,1}=\mathrm{Id},\quad\Phi_{1,2} \Phi_{2,1}=e^{\sqrt{1/2}\varphi},\quad  \Phi_{1,3} \Phi_{3,1}=\mathrm{Id}+2\epsilon e^{\sqrt{2}\varphi}, \\ \Phi_{1,4} \Phi_{4,1}=e^{-\sqrt{1/2}\varphi}+2\sqrt{2}\epsilon \partial\varphi e^{\sqrt{1/2}\varphi}+2\epsilon^2e^{3\sqrt{1/2}\varphi}.
\end{align*}
\end{Remark}

%

\section{Combinatorics}
\label{Seq:Comb} Recall, that $\calL_{h,k}$ denotes the irreducible 	highest weight representation of $\Asl(2)$. In this section we will consider the representations $\calL_{l,k}$ where $k \in \mathbb{N}, l\in \mathbb{Z}$ and $0 \leq l \leq k$. These representations are integrable. 

We call the function $f\colon \mathbb{Z}\rightarrow \mathbb{Z}_{\geq 0}$ a $(l,k)$ configuration if 
\begin{enumerate}
\item $f(m)+f(m+1)\leq k$ \label{con:1}
\item  $f(2m+1)=k-l$, $f(2m)=l$, for $m<<0$ \label{con:2}
\item  $f(m)=0$, for $m>>0$ \label{con:3}
\end{enumerate}
The set of such configurations we denote by $\Sigma_{l,k}$. By $f_n$ we denote so called extremal configurations 
\[f_n(m)=\begin{cases} 0\qquad &\text{if $m>n$;} \\ l\qquad &\text{if $m\leq n$, $m$ is even;}\\ k-l \qquad &\text{if $m\leq n$, $m$ is odd.}
\end{cases} \]
In the table below we represent configuration $f_{2n}$
\begin{center}
\begin{tikzpicture}[scale=0.9]
\foreach \i in {0,3,6}
{\node [anchor= south] at ($(-\i+1.5,0)$) {\small $l$};
\node [anchor=south] at ($(-\i,0)$) {\small $k-l$};}
\foreach \i in {0,1.5,3}
\node [anchor=south] at ($(\i+3,0)$) {\small $0$};
\node [anchor=south] at ($(7.5,0)$) {$\cdots$};
\node [anchor=south] at ($(-7.5,0)$) {$\cdots$};
\draw (-8,0) -- (8,0);
\foreach \i in {-3,-2,-1,+1,+2}
\node [anchor=north] at ($(1.5*\i+1.5,0)$) {\small $2n\i$};
\node [anchor=north] at ($(1.5,0)$) {\small $2n$};
\node [anchor=north] at ($(6,0)$) {$\cdots$};
\node [anchor=north] at ($(-4.5,0)$) {$\cdots$};
\node [anchor=south] at ($(-9,0)$) {$f_{2n}(m)$};
\node [anchor=north] at ($(-9,0)$) {$m$};
\draw (-8.15,-0.6) -- (-8.15,0.6);
\draw (-9.7,0) -- (-8.3,0);
\end{tikzpicture}
\end{center}

By $v$ we denote the highest weight vector of $\calL_{l,k}$. Extremal vectors $v_n \in \calL_{l,k}$, $n \in \mathbb{Z}$ defined by the relations 
\begin{equation}\label{eq:vn=vn}
 v_0=v, \qquad v_{2n}=(e_{-2n})^lv_{2n-1}, \quad v_{2n-1}=(e_{-2n+1})^{k-l}v_{2n-2},
\end{equation}
The Weyl group of $\Asl(2)$ acts on $\calL_{l,k}$ and the set of vectors $\{v_n\}$ is an orbit of the highest weight vector $v$ under the action of this group. We formally define $v_{-\infty}$ such that 
\begin{align*}
v_{2n}&=(e_{-2n})^l(e_{-2n+1})^{k-l}(e_{-2n+2})^l(e_{-2n+3})^{k-l}\cdots (e_{-2m})^l(e_{-2m+1})^{k-l}\cdots v_{-\infty}\\
v_{2n-1}&=(e_{-2n+1})^{k-l}(e_{-2n+2})^l(e_{-2n+3})^{k-l}(e_{-2n+4})^l\cdots (e_{-2m})^l(e_{-2m+1})^{k-l}\cdots v_{-\infty}
\end{align*}
Clearly these formulas agree with the equations \eqref{eq:vn=vn}. 
Due to condition \eqref{con:1} for any $f \in \Sigma_{l,k}$ there exist $n \in \mathbb{Z}$ such that for any $m >n$, $f(-2m)=l, f(-2m-1)=k-l$. Then $f$ differs from $f_n$ only in finite number of $m$ and we can define the semiinfinite product $\prod (e_{-m})^{f(m)} v_\infty$ by use of action of finite product of $e_m$ on $v_n$.
\begin{Theorem}[{\cite[Prop 2.6.1']{Feigin:1993qr}}]
The vectors of the form $\prod e_{-m}^{f(m)} v_\infty$ form a basis in the space $\calL_{l,k}$.
\end{Theorem}
As a consequence we can find the character of $\calL_{l,k}$. For any $f \in \Sigma_{l,k}$ we define the $q$-weight
\[w_q(f)=- \sum_{m< 0} (2m+1)(k-l-f(2m+1)) - \sum_{m< 0} 2m (l-f(2m))+ \sum_{m\geq 0} mf(m)\]
Due to conditions \eqref{con:2} and \eqref{con:3} this sum is finite. Clearly it is just the difference between $L_0$ gradings of the $\prod (e_{-n})^{f(n)} v_\infty$ and $v_0$. Since $L_0v=\frac{l(l+2)}{4(k+2)}v$ we have
\[ \chi(\calL_{l,k})=q^{\frac{l(l+2)}{4(k+2)}}\sum_{f \in \Sigma_{l,k}}q^{w_q(f)}.\]

Now we come to the main point of the section. We decompose the set $\Sigma_{l,k}$ as $\Sigma_{l,k}=\sqcup \Sigma_{l,k}^r$, where $\Sigma_{l,k}^r$ consists of $(l,k)$ configurations such that $f(0)=r$. It is clear that $\Sigma_{l,k}^r=\Sigma_{k}^{+, k-r}\times \Sigma_{l,k}^{-, k-r}$, where $\Sigma_{k}^{+,k-r}$ consists of functions $f\colon \mathbb{N} \rightarrow \mathbb{Z}_{\geq 0}$ such that $f(1)\leq k-r$ and conditions \eqref{con:1} and \eqref{con:3} hold and $\Sigma_{l,k}^{-,k-r}$ consists of functions $f\colon -\mathbb{N} \rightarrow \mathbb{Z}_{\geq 0}$ such that $f(-1)\leq k-r$ and conditions \eqref{con:1} and \eqref{con:2} hold. Therefore
\begin{equation}\label{eq:chi=chichi}
q^{-\frac{l(l+2)}{4(k+2)}}\cdot \chi(\calL_{l,k})=\sum_{f \in \Sigma_{l,k}}q^{w_q(f)}= \sum_{0\leq r \leq k} \left(\sum_{f \in \Sigma_{k}^{+, k-r}}q^{w_q(f)}\right)\cdot \left(\sum_{f \in \Sigma_{l,k}^{-, k-r}}q^{w_q(f)}\right) 
\end{equation}

It was proven in \cite[Prop. 5]{Feigin:1993} that the characters of the $(2,2k+3)$ minimal model representations have the form
\begin{equation}\label{eq:chi22r+1}
\chi^{2/(2k+3)}_{1,r}=q^{\Delta(P_{1,r},b_{2/(2k+3)})} \!\!\!\! \sum_{f \in \Sigma_{k}^{+,r-1}}\!\!\!q^{w_q(f)}=q^{\Delta(P_{1,r},b_{2/(2k+3)})} \!\!\!\!\!\!\sum_{n_1,n_2,\dots,n_{k}}\!\!\!\frac{q^{\sum_{i,j=1}^k n_in_j \mathrm{min}(i,j)+\sum_{j=r}^k (j-r+1)n_j}}{(q)_{n_1}\ldots (q)_{n_k}},
\end{equation}
where $b_{p/p'}^2=-p/p'$ and $(q)_n=\prod_{j=1}^n(1-q^j)$. This is the algebraic meaning of the first multipliers in  \eqref{eq:chi=chichi}. 

For the second multipliers let us consider $k=1$. In this case sums over $\Sigma_{l,k}^{-,r}$ can be simply rewritten in ``fermionic'' form. Due to e.g. \cite[eq. (8) and Theorem 2.3]{Feigin:1995} these fermionic formulas equal to the  characters of the representations of $(3,5)$ minimal model
\begin{align*}
\chi^{3/5}_{1,1}&=q^{\Delta(P_{1,1},b_{3/5})}\sum_{n=0}^\infty\frac{q^{n^2+n}}{(q)_{2n}}=q^{\Delta(P_{1,1},b_{3/5})}\sum_{f \in \Sigma_{1,1}^{-,0}}q^{w_q(f)},\\
\chi^{3/5}_{1,2}&=q^{\Delta(P_{1,2},b_{3/5})}\sum_{n=0}^\infty\frac{q^{n^2}}{(q)_{2n}}=q^{\Delta(P_{1,2},b_{3/5})}\sum_{f \in \Sigma_{0,1}^{-,1}}q^{w_q(f)},\\
\chi^{3/5}_{1,3}&=q^{\Delta(P_{1,3},b_{3/5})}\sum_{n=0}^\infty\frac{q^{n^2+n}}{(q)_{2n+1}}=q^{\Delta(P_{1,3},b_{3/5})}\sum_{f \in \Sigma_{1,1}^{-,1}}q^{w_q(f)},\\
\chi^{3/5}_{1,4}&=q^{\Delta(P_{1,4},b_{3/5})}\sum_{n=0}^\infty\frac{q^{n^2+2n}}{(q)_{2n+1}}=q^{\Delta(P_{1,4},b_{3/5})}\sum_{f \in \Sigma_{0,1}^{-,0}}q^{w_q(f)}.
\end{align*}
Therefore, for $k=1$ first multipliers in \eqref{eq:chi=chichi} are equal to the characters of $(2,5)$ minimal model (due to \eqref{eq:chi22r+1}) and the second equal to characters  of $(3,5)$ minimal model. So for $k=1$ the combinatorial identity \eqref{eq:chi=chichi} is equivalent to \eqref{eq:chi_sl2Vir2535}.

For $k>1$ the sums over $\Sigma_{l,k}^{-,r}$ should coincide with certain coset characters. It is natural to conjecture the isomorphism of product these coset algebra
with minimal model $\calM_{2/(2k+3)}$ and Urod algebra for $\Asl(2)$ on the level $k$ (similarly to the Theorem \ref{th:Ur}).

\section{Functional equations} \label{Se:Equations}
In this section we use the results of the Section \ref{Se:Res} for the functional equations on the conformal blocks. Then we explain the relation between these equations and and Nakajima-Yoshioka blow--up equations mentioned in the beginning of the paper \eqref{eq:blowr2}.

\subsection{Whittaker vector.} First we need to recall the definition of the Whittaker (or Gaiotto) limit of conformal block.

For the Verma module $\rmV_{\textsc{p},b}$ the Whittaker vector $\mathsf{W}_{\textsc{p},b}=\sum_{N=0} w_{\textsc{p},b,N}q^{N/2}$, where $w_{\textsc{p},b,N}\in \rmV_{\textsc{p},b}$, $L_0 w_{\textsc{p},b,N}=(\Delta+N)w_{\textsc{p},b,N}$ defined by the equations
\[ L_1 w_{\textsc{p},b,N}=w_{\textsc{p},b,N-1}, \quad \quad L_2 w_{\textsc{p},b,N}=0.\]
These equations can be simply rewritten as $L_1 \mathsf{W}_{\textsc{p},b}=q^{1/2}\mathsf{W}_{\textsc{p},b}$, $L_2 \mathsf{W}_{\textsc{p},b}=0$. \footnote{More general Whittaker vectors defined by the relations $L_1\mathsf{W}=\alpha\mathsf{W}$, $L_2\mathsf{W}=\beta\mathsf{W}$, for generic $\alpha$, $\beta$. The Whittaker vector for $\beta=0$ used in our paper sometimes called the Gaiotto vector }
It is easy to see that for generic $P,b$ the Whittaker vector $\mathsf{W}_{\textsc{p},b}$ exists and unique up to normalization. We will always use normalization of $\mathsf{W}_{\textsc{p},b}$ such that $\langle w_{\textsc{p},b,0},w_{\textsc{p},b,0}\rangle=1$, where $\langle\cdot,\cdot\rangle$ is a Shapovalov form in Verma module $\rmV_{\textsc{p},b}$.

The Whittaker limit of the four point conformal block is defined by
\begin{align}
\mathbb{F}(P,b;q)=\langle \mathsf{W}_{\textsc{p},b},\mathsf{W}_{\textsc{p},b}\rangle = \sum_{N=0}^{\infty}  \langle w_{\textsc{p},b,N},w_{\textsc{p},b,N}\rangle q^N \label{eq:WhitF}
\end{align}

Now we consider the representation $U_1 \otimes \bbL_{\textsc{p},b}$ of the algebra $\calA_b=\calU\otimes \calM_b$. Consider the vector $v_{\sqrt{1/2}}\otimes \mathsf{W}_{\textsc{p},b}(q)\in U_1 \otimes \bbL_{\textsc{p},b}$ (recall that after the shift of the grading $v_{\sqrt{1/2}}$ become the highest vector of $U_1$). It follows from the formulas \eqref{eq:Tb1},\eqref{eq:Tb2} that
\begin{align*}
L_1^{b_1}\left(v_{\sqrt{1/2}}\otimes \mathsf{W}_{\textsc{p},b}(q)\right)=q^{1/2}\frac{b^{-1}}{b^{-1}-b} \left( v_{\sqrt{1/2}}\otimes \mathsf{W}_{\textsc{p},b}(q)\right),  \qquad  L_2^{b_1}\left(v_{\sqrt{1/2}}\otimes \mathsf{W}_{\textsc{p},b}(q)\right)=0 \\
L_1^{b_2}\left(v_{\sqrt{1/2}}\otimes \mathsf{W}_{\textsc{p},b}(q)\right)=q^{1/2}\frac{b}{b-b^{-1}} \left( v_{\sqrt{1/2}}\otimes \mathsf{W}_{\textsc{p},b}(q)\right),  \qquad  L_2^{b_2}\left(v_{\sqrt{1/2}}\otimes \mathsf{W}_{\textsc{p},b}(q)\right)=0. 
\end{align*}
Therefore, using the decomposition of $U_1\otimes \bbL_{\textsc{p},b}$ from Theorem~\ref{th:Rep} we have the decomposition
\begin{equation}\label{eq:w=ww}
v_{\sqrt{1/2}}\otimes \mathsf{W}_{\textsc{p},b}(q)=\sum_{k \in \mathbb{Z}} \frac{q^{k^2/2}}{\sqrt{l_k(P,b)}} \left(\mathsf{W}_{P_1+kb_1, b_1}\left(\beta_1q\right)\otimes \mathsf{W}_{P_2+kb_2^{-1}, b_2}\left(\beta_2q\right)\right),
\end{equation}
where $\beta_1=\dfrac{b^{-2}}{(b^{-1}-b)^2}$,  $\beta_2=\dfrac{b^{2}} {(b-b^{-1})^2},$ the degrees $k^2/2$ defined by the difference in $L_0$ grading (see \eqref{eq:Pb}) \[k^2=\Delta(P_1+kb_1,b_1) +\Delta(P_2+kb_2^{-1},b_2) -\Delta(P_1,b_1) +\Delta(P_2,b_2)\] and $l_k(P,b)$ are unknown coefficient. Taking the norm of the right and left  sides of \eqref{eq:w=ww} we get
\begin{equation}\label{eq:F=FF}
\mathbb{F}(P,b;q)=\sum_{k \in \mathbb{Z}} \frac{q^{k^2}}{l_k(P,b)} \cdot \mathbb{F}\left(P_1+kb_1,b_1;\beta_1q\right)\cdot \mathbb{F}\left(P_2+kb_2^{-1},b_2;\beta_2q\right)
\end{equation}
The factors $l_k(P,b)$ in principle are determined by the equation \eqref{eq:w=ww}. Clearly $l_0(P,b)=1$. See also discussion in the end of this section.

\subsection{Differential equations.}
Let us consider the operator $H=bL_0^{b_1}+b^{-1}L_0^{b_2}$. It follows from the formulas \eqref{eq:Tb1},\eqref{eq:Tb2} that the corresponding local operator have the form
\begin{multline}\label{eq:H(z)}
bT_{b_1}+b^{-1}T_{b_2}=\frac{b+b^{-1}}{2\epsilon}e^{-\sqrt{2}\varphi} +\frac{b+b^{-1}}{2}(\partial \varphi)^2 +(b+b^{-1})\epsilon (\partial\varphi)^2e^{\sqrt{2}\varphi} \\ -\frac{2\epsilon^2}{b+b^{-1}}e^{2\sqrt{2}\varphi}  -\frac{2\epsilon}{b+b^{-1}} T_{b}e^{\sqrt{2}\varphi} 
\end{multline}
Define the function $\widehat{\mathbb{F}}$ by
\begin{align*}
\widehat{\bbF}(P,b;q,t)=\sum_{k=0}^{\infty} \widehat{\bbF}_m (P,b;q) \frac{t^m}{m!}=\left\langle v_{\sqrt{1/2}}\otimes \mathsf{W}_{\textsc{p},b} ,e^{tH} \left( v_{\sqrt{1/2}}\otimes \mathsf{W}_{\textsc{p},b}\right) \right\rangle
\end{align*}
It is clear from the definition of the operator $H$ that
\begin{multline*}
\widehat{\bbF}(P,b;q,t)=\sum_{k \in \mathbb{Z}} \frac{q^{k^2} }{l_k(P,b)} \cdot e^{t 	b\Delta^1_k }\mathbb{F}\left(P_1+kb_1,b_1;\beta_1qe^{tb} \right) \cdot e^{tb^{-1}\Delta^2_k}\mathbb{F}\left(P_2+kb_2^{-1},b_2;\beta_2q e^{tb^{-1}}\right),
\end{multline*}
where $\Delta^1_k=\Delta(P_1+kb_1,b_1)$ and $\Delta^2_k=\Delta(P_2+kb_2^{-1},b_2)$. Clearly $\widehat{\bbF}_0(P,b;q)=\bbF(P,b;q)$. In~order to write analogues formulas for $\widehat{\bbF}_m (P,b;q)$ we will use generalized Hirota-differential~\cite{NakajimaYoshioka}
\begin{equation*}
   \left(D^{(\epsilon_1,\epsilon_2)}_x\right)^m (f\cdot g)
   = \left.(\frac{d}{dy})^m f(x + \epsilon_1 y) g(x + \epsilon_2 y)\right|_{y = 0}
\end{equation*}
Therefore
\begin{equation}\label{eq:hatF}
\widehat{\bbF}_m (P,b;q)=\sum_{k \in \mathbb{Z}} \frac{q^{1/4-\Delta(P,b)}}{l_k(P,b)}\left( D^{(b,b^{-1})}_{\log q}\right)^m \left(q^{\Delta_k^1} \mathbb{F}\left( P_1+kb_1,b_1;\beta_1q\right)\cdot q^{\Delta_k^2} \mathbb{F}\left(P_2+kb_2^{-1},b_2;\beta_2q\right)\right), 
\end{equation}
where we used that $\Delta_k^1+\Delta_k^2=\Delta(P,b)-1/4+k^2$. 
		
On the other hand we can use the formula $\eqref{eq:H(z)}$ and apply $H$ to $v_{\sqrt{1/2}}\otimes \mathsf{W}_{\textsc{p},b}$. It is easy to see that
\[H \left(v_{\sqrt{1/2}}\otimes \mathsf{W}_{\textsc{p},b} \right)= \frac{b+b^{-1}}4\left(v_{\sqrt{1/2}}\otimes \mathsf{W}_{\textsc{p},b}\right)+\frac{-2\epsilon q^{1/2}}{b+b^{-1}} \left(v_{3/\sqrt{2}}\otimes \mathsf{W}_{\textsc{p},b}\right)\]
But the vectors $v_{3/\sqrt{2}}$ and $v_{1/\sqrt{2}}$ are orthogonal, hence we have
\begin{equation}
\widehat{\bbF}_1 (P,b;q)=\left\langle v_{\sqrt{1/2}}\otimes \mathsf{W}_{\textsc{p},b} ,H \left( v_{\sqrt{1/2}}\otimes \mathsf{W}_{\textsc{p},b}\right) \right\rangle=\frac{b+b^{-1}}4 {\bbF} (P,b;q) \label{eq:F1}
\end{equation}
Similarly applying $H$ one can prove that
\begin{align}
&\widehat{\bbF}_2 (P,b;q)=\left(\frac{b+b^{-1}}4\right)^2 {\bbF} (P,b;q), \qquad
\widehat{\bbF}_3 (P,b;q)=\left(\frac{b+b^{-1}}4\right)^3 \bbF (P,b;q) \label{eq:F23}\\
&\widehat{\bbF}_4 (P,b;q)=\left(\left(\frac{b+b^{-1}}4\right)^4 -2 q\right) \bbF (P,b;q) \label{eq:F4}\\
&\widehat{\bbF}_5 (P,b;q)=\left(\left(\frac{b+b^{-1}}4\right)^5 -\frac{17}{2}(b+b^{-1}) q\right) \bbF (P,b;q) \label{eq:F5} \\
&\widehat{\bbF}_6 (P,b;q)=\left(\left(\frac{b+b^{-1}}4\right)^6 -\frac{183(b+b^{-1})^2}{8} q \right){\bbF} (P,b;q)+ 8 q^{3-\Delta(P,b)} \partial_q \left(q^{\Delta(P,b)} \bbF (P,b;q)\right) \label{eq:F6}
\end{align}

One can easily see the general statement for the structure of the $\widehat{\bbF}_m$

\begin{Proposition} For any $m \geq 0$, the function $\widehat{\bbF}_m$ is a linear combination of the derivatives of $q^{-\Delta(P,b)} (\partial_q)^l \left(q^{\Delta(P,b)} \bbF (P,b;q)\right)$ with coefficients which are polynomials in $q$ and $b+b^{-1}$.
\end{Proposition}

This fact follows from the following two observations. First, in the  $\rmV_{\textsc{p},b}$ module the scalar product $\left\langle \mathsf{W}_{\textsc{p},b} ,L_{a_1}L_{a_2}\cdots L_{a_k} \mathsf{W}_{\textsc{p},b} \right\rangle$ where $a_1\leq a_2\leq \dots \leq a_k$ vanishes if $a_1<-1$ or $a_k>1$. And the nonzero products equals
\[ \left\langle \mathsf{W}_{\textsc{p},b} ,L_{-1}^{l}L_{0}^{l}L_{1'}^{l''} \mathsf{W}_{\textsc{p},b} \right\rangle = q^{(l+l'')/2}q^{-\Delta(P,b)+l'} (\partial_q)^{l'} \left(q^{\Delta(P,b)} \bbF (P,b;q)\right).\]

Second, we will have no $(b+b^{-1})$ in the denominator. The only possible origin of such denominators are the operators $-\frac{2\epsilon^2}{b+b^{-1}}e^{2\sqrt{2}\varphi}$ and  $-\frac{2\epsilon}{b+b^{-1}} T_{b}e^{\sqrt{2}\varphi}$. But this operators increase $\varphi_0$ grading. Therefore in order to have non zero scalar product the actions of such operators should be accompanied by we the action of the operator $\frac{b+b^{-1}}{2\epsilon}e^{-\sqrt{2}\varphi}$, which cancels the denominator.

\subsection{Geometric interpretation.} Recall the notation from the beginning of the introduction. The $\mathit{M}(2,N)$ denotes the moduli space of framed torsion free sheaves of rank 2 on $\mathbb{CP}^2$, $Z(\epsilon_1,\epsilon_2,a;q)$ denotes the generating function of the equivariant volumes of $\mathit{M}(2,N)$, $\epsilon_1,\epsilon_2,a_1,a_2$ are equivariant parameters, $a=(a_1-a_2)/2$. The AGT relation for Whittaker limit reads
\begin{equation}\label{eq:F=Z}
 \mathbb{F}(\frac{a}{\sqrt{\epsilon_1\epsilon_2}},\sqrt{\frac{\epsilon_1}{\epsilon_2}};\frac{q}{(\epsilon_1\epsilon_2)^2})= Z(\epsilon_1,\epsilon_2,a;q),
\end{equation}
This relation is proven now (e.g. follows from the more general results proven in \cite{Fateev:2009}, \cite{Alba:2010},\cite{Schiffmann:2012},\cite{Maulik:2012}). Using this identity we see that equation \eqref{eq:F=FF} coincides with the Nakajima-Yoshioka blow--up equation \eqref{eq:blowr2} up to factors $l_k$ and $\mathrm{l}_k$, which we comment below.

The insertion of the operator $H$ geometrically equivalent to multiplication of the integrand by the $\mu(C)+(b+b^{-1})/4$, where $C \subset \hat{\mathbb{CP}^2}$ is an exceptional divisor. The $\mu(C)$ is a cohomology class on $\hat{M}$ defined in \cite[p. 22]{NakajimaYoshioka}, the dual homology cycle consists of bundles that restrict to $C$ in a non-trivial way. 

We rewrite \eqref{eq:F1} and \eqref{eq:F23} substituting $\widehat{\bbF}_m$ and $\bbF=\widehat{\bbF}_0$ from equation \eqref{eq:hatF}. We obtain bilinear equations on functions $\mathbb{F}\left( P_1+kb_1,b_1;\beta_1q\right)$ and $\mathbb{F}\left(P_2+kb_2^{-1},b_2;\beta_2q\right)$. After AGT substitution \eqref{eq:F=Z}, these equations have the same form as the equation (6.14) in \cite{NakajimaYoshioka} for the $r=2$. 

On the other hand, one can consider these bilinear equation on $\bbF$ as the linear equations on $l_k$. These equations determines $l_k$ and therefore we have $l_k=\mathrm{l}_k$. The factors $\mathrm{l}_k$ in Nakajima-Yoshioka equations were determined geometrically in \cite{NakajimaYoshioka} and have the form
\begin{equation} \label{eq:ln}
l_k(P,b)=\mathrm{l}_k(P,b)=\prod_{i,j \geq 0, \, i+j< 2k} (-2P-i b- j b^{-1}) \prod_{i,j \geq 1, \, i+j< 2k} (2P+i b+ j b^{-1}).
\end{equation}

\begin{Remark}\label{Re:lk}
In the discussion above we used both AGT relation and Nakajima-Yoshioka equations in order to get \eqref{eq:ln}. In looks like the coefficients $l_k(P,b)$ can be determined by methods of conformal field theory as in \cite[Sec. 3.3]{Bershtein:2015}. These would give an independent proof of 
AGT relation or Nakajima-Yoshioka equations.
\end{Remark}

The bilinear equations obtained from \eqref{eq:F4},\eqref{eq:F4},\eqref{eq:F6} should be equaivalent to higher Nakajima-Yoshioka equations, which follow from the results of the paper \cite{NakajimaYoshioka5}.

\section{Conclusion}
In this paper we considered the very concrete vertex algebras, with explicit formulas. Possible generalizations may provide some understanding of the subject.

\subsection{}One can ask whether the product $\calM_{b_1}\otimes \calM_{b_2}$ can be extended by fields $\Phi_{1,3}\otimes \Phi_{3,1}$ such that the result is well defined vertex algebra. It can be argued that this can be done if the central charges are related by the equation $b_1^2+b_2^{-2}=n$, where $n \in \mathbb{Z}$. By AGT duality the $|n|>1$ case corresponds to the instanton counting on the Hirzebruch surface \cite{Bruzzo:2011}. 

More general vertex algebras can be constructed as an extension of the product of $n$ algebras 
$\calM_{b_1}\otimes \calM_{b_2} \otimes \dots \otimes \calM_{b_p}$, where central charges related by the relation $b_i^2+b_{i+1}^{-2}=n_i$, where $n_i\in \mathbb{Z}$.

Let us mention two special cases. If $n=-2$ then the extended product $\calM_{b_1}\otimes \calM_{b_2}$ corresponds to the instanton counting on the resolution of $\mathbb{C}^2/\mathbb{Z}_2$. This vertex algebra  have an isomorphic description similar  Theorem \ref{th:Th1} b). Namely the product $\calM_{b_1}\otimes \calM_{b_2}$ extended by fields $\Phi_{1,2}\otimes \Phi_{2,1}$ is isomorphic to the product $\calF\otimes\mathcal{SM}_c$, where $\mathcal{SM}_c$ is a Super-Virasoro vertex superalgebra and $\mathcal{F}$ is the Majorana fermion superalgebra (see \cite{BBFLT}, and references therein).

If $n=0$ then the corresponding central charges are related by $c_1+c_2=26$. In this case one multiply by the ghost representation $\Lambda^{bc}$ and compute BRST cohomology. These cohomology classes are physical states in the chiral part of the Liouville gravity.

\subsection{} 
The Urod algebra $\calU$ considered in this paper is a deformation of the $\Asl(2)_1$ as a chiral algebra. It is natural to expect the existence of the such Urod deformations of the $\Asl(r)_k$ for any $k$ and $r \in \mathbb{N}$. Clearly the Section \ref{Se:proof} constructions (combination of coset construction and Drinfeld--Sokolov reduction) have such generalization.

The Urod algebra for $\Asl(r)_1$ by AGT relation corresponds to the the blow--up equations for the $U(r)$ instantons. Similar to $r=2$ case this algebra have a subalgebra isomorphic to the product of two $W_{r}$ algebra corresponding to minimal models $(r,2r+1)$ and $(2r+1,r+1)$. This implies the character identities, similar to \eqref{eq:chi_sl2Vir2535}. The corresponding identity for $r=3$ and vacuum representation reads
\begin{multline*}
\chi(\calL_{0,0,1})=q^{-1}\left(\chi^{3/7}_{(0,0)}\cdot \chi^{7/4}_{(0,0)}+\chi^{3/7}_{(0,3\Lambda_1)}\cdot \chi^{7/4}_{(3\Lambda_1,0)}+\chi^{3/7}_{(0,3\Lambda_2)}\cdot \chi^{7/4}_{(3\Lambda_2,0)}+\right. 
\\ 
\left.+ \chi^{3/7}_{(0,\Lambda_1+\Lambda_2)}\cdot \chi^{7/4}_{(\Lambda_1+\Lambda_2,0)}+\chi^{3/7}_{(0,2\Lambda_1+2\Lambda_2)}\cdot \chi^{7/4}_{(2\Lambda_1+2\Lambda_2,0)}\right),
\end{multline*}
where $\calL_{0,0,1}$ is a vacuum level 1 representation of $\Asl(3)$, $\Lambda_1,\Lambda_2$ are fundamental $\mathfrak{sl}(3)$ weights.

As another particular case we discuss the Urod algebra for $\Asl(2)_k$, where $k \in \mathbb{N}$. From the section \ref{Se:proof} constructions and combinatorial arguments from the section \ref{Seq:Comb} follow that the $\Asl(2)_k$ Urod algebra contains a subalgebra which is isomorphic to the product of $(2,2k+1)$ Virasoro minimal model and coset minimal model. For example $\Asl(2)_2$ Urod algebra has a subalgebra which is the product of $(2,7)$ Virasoro algebra and $(3,7)$ Super Virasoro algebra. The geometrical meaning of these Urod algebras is unknown.


\section{Acknowledgments}
We thank  A. Belavin and H. Nakajima for interest to our work and discussions. 

The financial support from the Government of the Russian Federation within the framework of the implementation of the 5-100 Programme Roadmap of the National Research University Higher School of Economics is acknowledged.
The research about the bilinear equation was performed under a grant funded by Russian Science Foundation (project No. 14-12-01383)


\begin{thebibliography}{99}
\bibitem{Alba:2010}V.~A. Alba, V.~A. Fateev, A.~V. Litvinov, G.~M. Tarnopolsky, {\it {On combinatorial expansion of the conformal blocks arising from AGT conjecture}},  Lett. Math. Phys. {\bf 98} (2011) 33-64,
  [\href{http://arxiv.org/abs/1012.1312}{{\tt arXiv:1012.1312}}].

\bibitem{AGT:2009} L.~F. Alday, D.~Gaiotto, Y.~Tachikawa, \emph{Liouville Correlation   Functions from Four-dimensional Gauge Theories},   Lett. Math. Phys.
  {\bf 91} (2010) 167-197,[\href{http://arxiv.org/abs/0906.3219}{{\tt
  arXiv:0906.3219}}].


\bibitem{BPZ} A.~Belavin, A.~Polyakov, A.~Zamolodchikov, {\it Infinite Conformal Symmetry in Two-Dimensional Quantum Field Theory.} Nucl. Phys. {\bf B241} (1984), 333.

\bibitem{BBFLT}  A.~Belavin, M.~Bershtein, B.~Feigin, A.~Litvinov,  G.~Tarnopolsky, {\it Instanton moduli spaces and bases in coset conformal field theory}. Comm. Math. Phys. \textbf{319 1}, 269-301 (2013), [\href{http://arxiv.org/abs/1111.2803}{{\tt  arXiv:1111.2803}}]

\bibitem{BGG} J.~Bernstein, I.~Gel'fand, S.~Gel'fand, \textit{Category of $\mathfrak{g}$-modules} Funkts. Anal. Prilozh., \textbf{10 2} (1976), 1–8.

\bibitem{Bershtein:2015} M. Bershtein , A. Shchechkin, {\it Bilinear Equations on Painlev\'e $\tau$ Functions from CFT}. Comm. Math. Phys. \textbf{339 3}, 1021-1061 (2015), [\href{http://arxiv.org/abs/1406.3008}{{\tt  arXiv:1406.3008}}]


\bibitem{Bruzzo:2011}
U.~{Bruzzo}, R.~{Poghossian}, and A.~{Tanzini}, {\it {Poincar{\'e} Polynomial
  of Moduli Spaces of Framed Sheaves on (Stacky) Hirzebruch Surfaces}},  {\em
  Commun. Math. Phys.} {\bf 304} (2011) 395--409,
  [\href{http://arxiv.org/abs/0909.1458}{{\tt arXiv:0909.1458}}].


\bibitem{Fateev:2009} V.A.~Fateev, A.V.~Litvinov, \textit{On AGT conjecture} JHEP \textbf{1002} (2010) 014 [\href{http://arxiv.org/abs/0912.0504}{{\tt arXiv:0912.0504}}].


\bibitem{Feigin:1995} B. Feigin, O. Foda, T. Welsh, \emph{Andrews-Gordon type identities from combinations of Virasoro characters}, Ramanujan J., \textbf{17 (1)}, (2008) 33-52 ; [\href{http://arxiv.org/abs/math-ph/0504014}{{\tt   arXiv:math-ph/0504014.}}]
.

\bibitem{Feigin:1990} B. Feigin, E. Frenkel, \emph{Quantization of the Drinfel′d-Sokolov reduction}. Phys. Lett. B, \textbf{246(1-2)} (1990) 75.

\bibitem{Feigin1:1990} B. Feigin, E. Frenkel, \emph{Affine Kac-Moody algebras, bosonization and resolutions}. Lett. Math. Phys. 19, (1990) 307-317.

\bibitem{Feigin:1993} B. Feigin, E. Frenkel, \emph{Coinvariants of nilpotent subalgebras of the Virasoro algebra and partition identities.} Adv. Sov. Math., \textbf{16}, (1993), 139–148 [\href{http://arxiv.org/abs/hep-th/9301039}{{\tt   arXiv:hep-th/9301039.}}]

\bibitem{Feigin Fuchs 1990} B.~Feigin, D.~Fuchs, {\it Representations of the Virasoro algebra.} Representations of Lie Groups and Related Topics, 465, Adv. Stud. Contemp. Math., 7, Gordon and Breach, New York, 1990.


\bibitem{Feigin:1993qr} B.L. Feigin, A.V. Stoyanovsky, Quasi-particles models for the representation of Lie algebras and geometry of flag
manifold, Funct. Anal. Appl. 28 (1994) 68–90, [\href{http://arxiv.org/abs/hep-th/9308079}{{\tt   arXiv:hep-th/9308079.}}]

\bibitem{FrenkelBook} E. Frenkel, D. Ben-Zvi, \emph{Vertex Algebras and Algebraic Curves}, Mathematical Surveys and Monographs 88,  American Mathematical Society 2004

\bibitem{Frenkel Kac} I. B. Frenkel, V. G. Kac, \emph{Basic representations of affine Lie algebras and dual resonance models}, Invent. Math. \textbf{62} (1980/81), no. 1.

\bibitem{ope.math} A. Fujitsu, \emph{ope.math: Operator product expansions in free field realizations of conformal field theory}, Comput. Phys. Commun. 79 (1994) 78-99.

\bibitem{Iohara Koga} K. Iohara, Y. Koga, {\it Representation theory of the Virasoro algebra}, Springer Monographs in Mathematics, London: Springer-Verlag London Ltd (2011)

\bibitem{KacKazhdan} V. G. Kac, D. A. Kazhdan, \textit{Structure of representations with highest weight of infinite-dimensional Lie algebras}, Adv. in Math. \textbf{34 1} (1979), 97–108.

\bibitem{KacWakimoto} V. G. Kac, M. Wakimoto, \textit{Modular invariant representations of infinite-dimensional Lie algebras and superalgebras} Proc. Natl. Acad. Sci. USA, \textbf{85} (1988), 4956- 
4960. 

\bibitem{GKO} P.~Goddard, A.~Kent, D.~Olive, \textit{Unitary representations of the Virasoro and super-Virasoro algebras} Comm. Math. Phys. \textbf{103 1} (1986), 105.

\bibitem{Maulik:2012} D.~Maulik, A.~Okounkov, \emph{Quantum Groups and Quantum Cohomology} [\href{http://arxiv.org/abs/1211.1287}{{\tt   arXiv:1211.1287}}]

\bibitem{NakajimaYoshioka} H.~Nakajima, K.~Yoshioka, \emph{Instanton counting on blowup. I. 4-dimensional pure gauge theory}, Inventiones mathematicae \textbf{162 2} (2005), 313-355 [\href{http://arxiv.org/abs/math/0306198}{{\tt
  arXiv:math/0306198}}].  

\bibitem{NakajimaYoshioka5} H.~Nakajima, K.~Yoshioka, \emph{Perverse coherent sheaves on blow-up. III. Blow-up formula from wall-crossing}, Kyoto J. Math. \textbf{51 2} (2011), 263 [\href{http://arxiv.org/abs/0911.1773}{{\tt
  arXiv:0911.1773}}].  


\bibitem{Schiffmann:2012} O. Schiffmann, E. Vasserot, \emph{Cherednik algebras, $W$ algebras and the equivariant cohomology of the moduli space of instantons on $\mathbb{A}^2$}, [\href{http://arxiv.org/abs/1202.2756}{{\tt
  arXiv:1202.2756}}].  
\end{thebibliography}
\end{document}